\RequirePackage{fix-cm}
\documentclass{svjour3}[2.0]

\usepackage[utf8]{inputenc}
\usepackage{units}
\usepackage{cite}
\usepackage{nicefrac}
\usepackage[numbers,sort&compress]{natbib}

\title{Fractional Calculus in Pharmacokinetics}%
\author{Pantelis~Sopasakis \and Haralambos~Sarimveis \and Panos~Macheras \and Aristides~Dokoumetzidis}
\institute{P. Sopasakis \at
	   KU Leuven, Department of Electrical Engineering (ESAT), STADIUS Center for Dynamical Systems, Signal Processing and Data Analytics, Kasteelpark Arenberg 10, 3001 Leuven, Belgium.\\
	   Tel.: +32-486-928034\\
	   \email{pantelis.sopasakis@kuleuven.be}
	   \and
	   H. Sarimveis \at
           School of Chemical Engineering, National Technical University of Athens, 9 Heroon Polytechneiou Street, 15780 Zografou Campus, Athens, Greece.\\
           Tel.: +30-210-7723237\\
           \email{hsarimv@central.ntua.gr}
           \and
	      P.~Macheras \at
              Department of Pharmacy,
	      University of Athens, Panepistimiopolis Zografou, 15784 Athens, Greece.\\
	      Tel: +30-210-7274026\\
	      \email{macheras@pharm.uoa.gr}
          \and
	      A.~Dokoumetzidis\at
              Department of Pharmacy,
	      University of Athens, Panepistimiopolis Zografou, 15784 Athens, Greece.\\
	      Tel: +30-210-7274122\\
	      \email{adokoum@pharm.uoa.gr}
}

 \journalname{J. Pharmacokinet. Pharmacodyn.}

\usepackage{amsmath,amssymb}
\usepackage{nicefrac}
\usepackage{graphics,graphicx}
\usepackage{geometry}
\usepackage{color}
\usepackage{hyperref}

\renewcommand{\Re}{\mathbb{R}}
\renewcommand{\d}{\mathrm{d}}
\newcommand{\D}{\mathrm{D}}
\newcommand{\N}{\mathbb{N}}
\newcommand{\E}{\mathcal{E}}
\renewcommand{\L}{\mathcal{L}}
\newcommand{\half}{\nicefrac{1}{2}}
\newcommand{\erf}{\operatorname{erf}}
\newcommand{\glD}{{{}^{\mathrm{gl}}\D}_t}
\newcommand{\glDD}{{^{\mathrm{gl}}\Delta}}
\newcommand{\rlI}{{{}_{0}^{\mathrm{rl}}\mathrm{I}}}
\begin{document}

\maketitle

\begin{abstract}
We are witnessing the birth of a new variety of pharmacokinetics 
where non-integer-order differential equations are employed to study the time course 
of drugs in the body: this is dubbed ``fractional pharmacokinetics''. 
The presence of fractional kinetics has important clinical implications such as the lack 
of a half-life, observed, for example with the drug amiodarone and the associated 
irregular accumulation patterns following constant and multiple-dose administration. 
Building models that accurately reflect this behaviour is essential for the design of less 
toxic and more effective drug administration protocols and devices. 
This article introduces the readers to the theory of fractional pharmacokinetics
and the research challenges that arise. After a short introduction 
to the concepts of fractional calculus, and the main applications that have appeared in 
literature up to date, we address two important aspects. First, numerical methods
that allow us to simulate fractional order systems accurately and second, optimal control 
methodologies that can be used to design dosing regimens to individuals and populations.
\end{abstract}


\section*{Introduction}

\subsection*{Background}
Diffusion is one of the main mechanisms of various transport processes in living 
species and plays an important role in the distribution of drugs in the body. 
Processes such as membrane permeation, dissolution of solids and dispersion in 
cellular matrices are considered to take place by diffusion. Diffusion is 
typically described by Fick's law, which, in terms of the pharmacokinetics of drugs,
gives rise to exponential washout curves that have a 
characteristic time scale, usually expressed as a half-life. 
However, in the last few decades, strong experimental evidence has suggested that 
this is not always true and diffusion processes that deviate from this law have 
been observed. These are either faster (super-diffusion) or slower (sub-diffusion) 
modes of diffusion compared to the regular case~\cite{art01,IONESCU2017141}. 

Such types of diffusion give rise to kinetics that are referred to as 
\textit{anomalous} to indicate the fact that they stray from the standard 
diffusion dynamics~\cite{art01}. Moreover, anomalous kinetics can also result 
from reaction-limited processes and long-time trapping. 
Anomalous kinetics introduces memory effects into the distribution process that need 
to be accounted for to correctly describe it. 
An distinctive feature of anomalous power-law kinetics is that it 
lacks a characteristic time scale contrary to exponential kinetics. A mathematical formulation that describes such anomalous kinetics is 
known as \textit{fractal kinetics}~\cite{kopelman,art05,Per10} where explicit power functions of time in the form of 
time-dependent coefficients are used to account for the memory effects replacing the 
original rate constants. In pharmacokinetics several datasets have been characterised by 
power laws, epmpirically~\cite{art02,art03}, while the first article that utilised 
fractal kinetics in pharmacokinetics was~\cite{art05}. 
Molecules that their  kinetics presents power law behaviour incude those distributed in deeper 
tissues, such as amiodarone~\cite{art04} and bone-seeking elements, such as calcium, lead, strontium and plutonium~\cite{art05,art06}. 

An alternative theory to describe anomalous kinetics employs fractional 
calculus~\cite{art07,Pod99}, which introduces derivatives and integrals of 
fractional order, such as half or three quarters. Although fractional calculus 
was introduced by Leibniz more than 300 years ago, it is only within the last 
couple of decades that real-life applications have been explored%
~\cite{art09,art10,art11}. It has been shown that differential equations with 
fractional derivatives (FDEs) describe experimental data of anomalous diffusion 
more accurately. Although fractional-order derivatives were first introduced as a novel mathematical concept with
unclear physical meaning, nowadays, a clear connection between diffusion over fractal spaces (such 
as networks of capillaries) and fractional-order dynamics has been established%
~\cite{BUTERA2014146,CHEN20101754,METZLER199413,COPOT20149277}. In particular, in anomalous diffusion
the standard assumption that the mean square displacement $\langle x^2 \rangle$
is proportional to time does not hold. Instead, it is \mbox{$\langle x^2 \rangle\sim t^{\nicefrac{2}{d}}$},
where $d\neq 2$ is the associated fractal dimension~\cite{Gmachowski2015,KlaSok2011,Eirich1990}.

Fractional kinetics~\cite{fractal2fractional} was introduced in the pharmaceutical literature in~\cite{art13} 
and the first example of fractional pharmacokinetics which appeared in that article 
was amiodarone, a drug known for its power-law kinetics~\cite{art03}. Since then, 
other applications of fractional pharmacokinetics have appeared in the literature. Kytariolos \textit{et al.}~\cite{KytMacDok10} presented an application of fractional kinetics for the development of nonlinear in vitro-in vivo correlations.
Popovi\'c \textit{et al.} have presented several applications of fractional 
pharmacokinetics to model drugs, namely for diclofenac~\cite{art20}, 
valproic acid~\cite{art21}, bumetanide~\cite{art22} and methotrexate~\cite{art23}. 
Copot \textit{et al.} have further used a fractional pharmacokinetics model for 
propofol~\cite{proc24}. In most of these cases the fractional model has been compared with 
an equivalent ordinary PK model and was found superior. 
FDEs have been proposed to describe drug response too, apart from their kinetics. 
Verotta has proposed several alternative fractional PKPD models that are capable of 
describing pharmacodynamic times series with favourable properties~\cite{art25}. Although 
these models are empirical, i.e., they have no mechanistic basis, they are attractive since 
the memory effects of FDEs can link smoothly the concentration to the response with a 
variable degree of influence, while the shape of the responses generated by fractional PKPD 
models can be very flexible and parsimonious (modelled using few parameters).

Applications of FDEs in pharmacokinetics fall in the scope of the newly-coined discipline 
of \textit{mathematical pharmacology}~\cite{vanderGraaf2015} to which this 
issue is dedicated. Mathematical pharmacology utilises applied mathematics, beyond the standard 
tools used commonly in pharmacometrics, to describe drug processes in the body and assist 
in controlling them.

This paper introduces the basic principles of fractional calculus 
and FDEs, reviews recent applications of fractional calculus in pharmacokinetics 
and discusses their clinical implications. 
Moreover, some aspects of drug dosage regimen optimisation based on control theory
are presented for the case of pharmacokinetic models that follow fractional kinetics. 
Indeed, in clinical pharmacokinetics and therapeutic drug monitoring, dose optimisation 
is carried out usually by utilising Bayesian methodology. Optimal control theory is a 
powerful approach that can be used to optimise drug administration, which can handle 
complicated constraints and has not been used extensively for this task. The case of 
fractional pharmacokinetic models is of particular interest for control theory and poses 
new challenges.

A problem that has hindered the applications of fractional calculus is the lack 
of efficient general-purpose numerical solvers, as opposed to ordinary differential equations. 
However, in the past few years significant progress has been achieved in this area. 
This important issue is discussed in the paper and some techniques for numerical solution of 
FDEs are presented.

\subsection*{Fractional derivatives}
Derivatives of integer order $n\in\N$, of functions $f:\Re\to\Re$ are well 
defined and their properties have been extensively studied in real analysis. 
The basis for the extension of such derivatives to real orders $\alpha\in\Re$
commences with the definition of an integral of order $\alpha$ which hinges 
on Cauchy's iterated ($n$-th order) integral formula and gives rise to 
the celebrated Riemann-Liouville integral $\rlI^{\alpha}_t$~\cite{Hennion2013}:
\begin{align}
\label{eq:riemann-liouville-integral}
\rlI^{\alpha}_t f(t) = {}_{0}\D_{t}^{-\alpha}f(t) 
   = \tfrac{1}{\Gamma(\alpha)}\int_0^t (t-\tau)^{\alpha-1}f(\tau)\d\tau,
\end{align}
where $\Gamma$ is the Euler gamma function. 
In Eq.~\eqref{eq:riemann-liouville-integral} we assume that $f$ is such that the involved
integral is well-defined and finite. This is the case if $f$ is continuous everywhere except  
at finitely many points and left-bounded, i.e., $\{f(x); 0 \leq x\leq t\}$ is bounded for every $t\geq 0$.
Note that for $\alpha\in\N$, $\rlI^{\alpha}_t f$ is equivalent to the ordinary $\alpha$-th order integral.
Hereafter, we focus on derivatives of order $\alpha \in [0,1]$ as only these
are of interest to date in the study of fractional pharmacokinetics.

The left-side subscript $0$ of the ${}_{0}\D_t^{\alpha}$ and $\rlI_t^{\alpha}$ operators, denotes the lower end of the integration 
limits, which in this case has been assumed to be zero. However, alternative lower bounds can be 
considered leading to different definitions of the fractional derivative with slightly 
different properties. When the lower bound is 
$-\infty$, the entire \textit{history} of the studied function is accounted for, which is considered preferable in some 
applications and is referred to as the \textit{Weyl} definition~\cite{art09}. 

It is intuitive to define a fractional-order
derivative of order $\alpha$ via ${}_{0}\D_{t}^{\alpha}= \D^1 {\,}_{0}\D_{t}^{\alpha-1} = \D^1 \,\rlI^{\alpha-1}$,
(where $\D^1$ is the ordinary derivative of order $1$), or equivalently
\begin{align}
{}_{0}\D_{t}^{\alpha}f(t) = \tfrac{\d}{\d t}\left[ 
  \tfrac{1}{\Gamma(1-\alpha)}\int_0^t \frac{f(\tau)}{(t-\tau)^{\alpha}}\d\tau
\right]
\end{align}
This is the \textit{Riemann-Liouville} definition --- one of the most popular 
constructs in fractional calculus.
One may observe that the fractional integration is basically a convolution between function $f$
and a power law of time, i.e., ${}_{0}\D_{t}^{-\alpha}f(t) = t^{\alpha-1}\ast f(t)$, where $\ast$ denotes
the convolution of two functions. This explicitly demonstrates the memory effects of the studied process. 
Fractional derivatives possess properties that are not straightforward or intuitive; 
for example, the half derivative of a constant $f(t) = \lambda$ with respect to $t$ does not vanish and instead is
$_{0}\D_t^{\nicefrac{1}{2}} \lambda = \nicefrac{\lambda}{\sqrt{\pi t}}$.

Perhaps the most notable shortcoming of the Riemann-Liouville definition with the $0$ lower bound is 
that when used in differential equations it gives rise to initial conditions that 
involve the fractional integral of the function and are difficult to interpret physically. 
This is one of the reasons the Weyl definition was introduced, but 
this definition may not be very practical for most applications either, as it involves an 
initial condition at time $-\infty$~\cite{art09,Sam+93}. 

A third definition of the fractional derivative, which 
is referred to as the \textit{Caputo} definition is preferable for most physical processes as 
it involves explicitly the initial condition at time zero. The definition is:
\begin{align}
{}_{0}^{\mathrm{c}}\D_{t}^{\alpha}f(t) 
  = \tfrac{1}{\Gamma(1-\alpha)}\int_0^t\frac{\tfrac{\d f(\tau)}{\d\tau}}{(t-\tau)^\alpha}\d\tau
\end{align}
where the upper-left superscript $\mathrm{c}$ stands for \textit{Caputo}. 
The Caputo derivative is interpreted as ${}_{0}^{\mathrm{c}}\D_{t}^{\alpha}= \rlI^{\alpha-1} \D^1$,
that is $\D^1$ and $\rlI^{\alpha-1}$ are composed in the opposite way compared
to the Riemann-Liouville definition. The Caputo derivative gives rise to initial 
value problems of the form
\begin{subequations}
\begin{align}
 {}_{0}^{\mathrm{c}}\D_{t}^{\alpha} x(t) &= F(x(t)),\\
 x(0) &= x_0.
\end{align}
\end{subequations}

This definition for the fractional derivative, apart from the more familiar 
initial conditions, comes with some properties similar to those of integer-order derivatives, 
one of them being that the Caputo derivative of a constant is in fact zero. 
Well-posedeness conditions, for the existence and uniqueness of solutions, 
for such fractional-order initial value problems are akin to those of integer-order problems~\cite{DENG20146},\cite[Chap.~3]{Pod99}.
The various definitions of the fractional derivative give different results, but these 
are not contradictory since they apply for different conditions and it is a 
matter of choosing the appropriate one for each specific application. All definitions 
collapse to the usual derivative for integer values of the order of differentiation. 

A fourth fractional derivative definition, which is of particular interest is the \textit{Gr\"{u}nwald-Letnikov},
which allows the approximate discretisation of fractional differential equations.
The Gr\"{u}nwald-Letnikov derivative, similar to integer order derivatives, is defined
via a limit of a \textit{fractional-order difference}. Let us first define
the Gr\"{u}nwald-Letnikov difference of order $\alpha$ and step-size $h$ of
a left-bounded function $f$ as
\begin{align}\label{eq:gl-diff}
\glDD^\alpha_h f(t) = \tfrac{1}{h^{\alpha}} \sum_{j=0}^{ \left\lfloor \nicefrac{t}{h}\right\rfloor}c_j^\alpha f(t-jh),
\end{align}
with $c^\alpha0=1$ and for $j\in\N_{>0}$
\begin{align}
c_j^\alpha=(-1)^{j}\prod_{i=0}^{j-1}\frac{\alpha-i}{i+1}.
\end{align}

The Gr\"{u}nwald-Letnikov difference operator leads to the definition of 
the \textit{Gr\"{u}nwald-Letnikov derivative} of order $\alpha$ which is 
defined as~\cite[Section~{20}]{Sam+93}
\begin{align}
\glD^\alpha f(t)=\lim_{h\to 0^+}\glDD^\alpha_h f(t),
\end{align}
insofar as the limit exists. This definition is of great importance in practice,
as it allows us to replace $\glD$ with $\glDD_h^\alpha$ provided that $h$ is 
adequately small, therefore it serves as an Euler-type discretisation of 
fractional-order continuous time dynamics. By doing so, we produce a discrete-time
yet infinite-dimensional approximation of the system since $\glDD_h^\alpha f(t)$ depends
on $f(t-jh)$ for all $j\in\N_{>0}$.
However, we may truncate this series up to a finite history $\nu$, defining 
\begin{align}
 \glDD_{h,\nu}^\alpha f(t) =\tfrac{1}{h^{\alpha}}
   \sum_{j=0}^{\min\{\nu,\left\lfloor \nicefrac{t}{h}\right\rfloor\}}
   c_j^\alpha f(t-jh).
\end{align}
This leads to discrete-time finite-memory approximations which are particularly 
useful as we shall discuss in what follows.

\section*{Fractional Kinetics}
The most common type of kinetics encountered in pharmaceutical literature are the so called ``zero order'' and ``first order''. 
Here the ``order'' refers to the order of linearity and is not to be confused with the  order of differentiation, i.e., a zero order process refers to a constant rate and a first order to a proportional rate. The fractional versions of these types of kinetics are presented below and take the form of fractional-order ordinary differential equations. Throughout this presentation the Caputo version of the fractional derivative is considered for reasons already explained.

\subsection*{Zero-order kinetics} 
The simplest kinetic model is the zero-order model where it is assumed that
the rate of change of the quantity $q$, expressed in [mass] units, is 
constant and equal to $k_0$, expressed in [mass]/[time] units. Zero-order 
systems are governed by differential equations of the form
\begin{align}
 \frac{\d q}{\d t} = k_0.
\end{align}
The solution of this equation with initial condition $q(0)=0$ is
\begin{align}
 q(t) = k_0 t.
\end{align}

The fractional-order counterpart of such zero-order kinetics can be obtained 
by replacing the derivative of order $1$ by a derivative of fractional order 
$\alpha$, that is
\begin{align}
 {}_{0}^{\mathrm{c}}\D_{t}^{\alpha}q(t) = k_{0,f},
\end{align}
where $k_{0,f}$ is a constant with units [mass]/[time]$^{\alpha}$. 
The solution of this equation for initial condition $q(0)=0$ is a power law
for $t>0$~\cite{Pod99}
\begin{align}
 q(t) = \frac{k_{0,f}}{\Gamma(\alpha+1)}t^\alpha.
\end{align}

\subsection*{First-order kinetics} 
The first-order differential equation, where the rate of change of 
the quantity $q$ is proportional to its current value, is described by
the ordinary differential equation
\begin{align}
 \frac{\d q}{\d t} = -k_1 q(t).
\end{align}
Its solution with initial condition of $q(0)=q_0$ 
is given by the classical equation of exponential decay
\begin{align}
 q(t) = q_0 \exp(-k_1 t).
\end{align}

Likewise, the fractional-order analogue of such first-order kinetics
is derived by replacing $\nicefrac{\d}{\d t}$ by the fractional-order
derivative ${}_{0}^{\mathrm{c}}\D_{t}^{\alpha}$ yielding the 
following fractional differential equation
\begin{align}
 \label{eq:fde-first-order}
 {}_{0}^{\mathrm{c}}\D_{t}^{\alpha}q(t) =- k_{1,f} q(t),
\end{align}
where $k_{1,f}$ is a constant with [time]$^{-\alpha}$ units. The solution 
of this equation can be found in most books or papers of the fast 
growing literature on fractional calculus~\cite{Pod99} and for initial condition 
$q(0)=q_0$ it is
\begin{align}
 \label{eq:fde-first-order-solution}
 q(t) = q_0 \E_\alpha(-k_{1,f}t^\alpha),
\end{align}
where $\E_\alpha$ is the Mittag-Leffler function~\cite{Pod99} which is defined as 
\begin{align}
 \E_\alpha(t) = \sum_{k=0}^{\infty}\frac{t^k}{\Gamma(\alpha k + 1)}.
\end{align}

The function $\E_\alpha(t)$ is a generalisation of the exponential function and it 
collapses to the exponential when $\alpha=1$, i.e., $\E_1(t)=\exp(t)$.
Alternatively, Eq.~\eqref{eq:fde-first-order-solution} can be reparametrised by introducing 
a time scale parameter with regular time units, $\tau_f = k_{1,f}^{-1/\alpha}$ and then 
becomes $q(t) = q_0 \E_\alpha(-(t/\tau_f)^\alpha)$. 
The solution of Eq.~\eqref{eq:fde-first-order} basically means that the fractional 
derivative of order $\alpha$ of the function $\E_\alpha(t^\alpha)$ is itself a 
function of the same form, exactly like the classic derivative of an exponential 
is also an exponential. 

\begin{figure}
 \centering
 \includegraphics[width=0.65\textwidth]{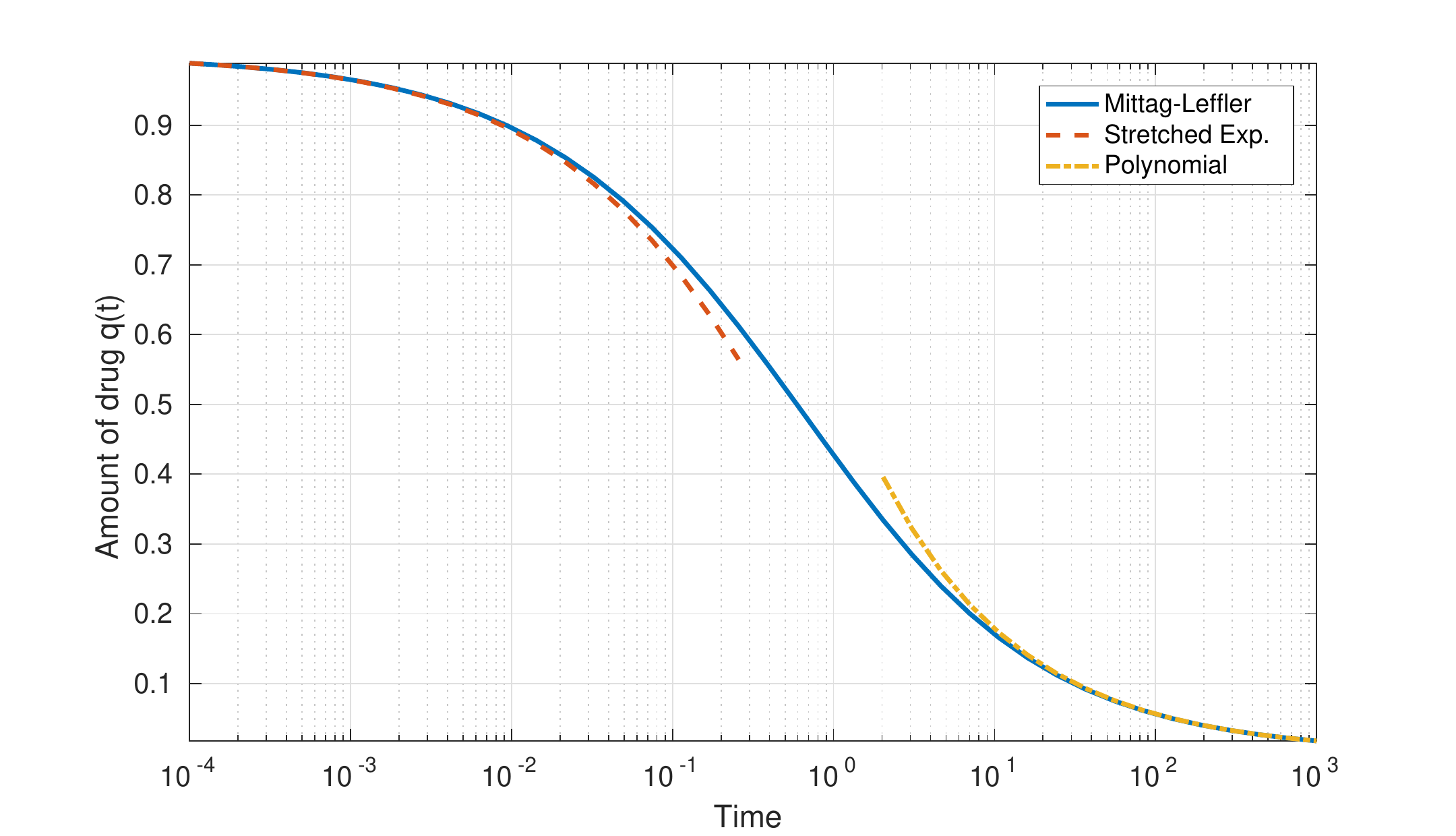}
 \caption{Amount-time profile according to Eq.~\eqref{eq:fde-first-order-solution} 
 for $\alpha=0.5$ (solid curve) along with its approximation at values of $t$ close to $0$
 by a stretched exponential function $\exp(-t^{\half}/\Gamma(1.5))$ (bottom dashed curve) and 
 an approximation at high values of $t$ by a power function $t^{-\half}/\Gamma(0.5)$ 
 (top dashed-dotted curve). Note that the time axis is logarithmic. The evolution of the 
 amount of drug starts as a stretched exponential and eventually shapes as a power function.}
 \label{fig:f1}
\end{figure}

It also makes sense to restrict $\alpha$ to values $0\leq \alpha \leq 1$, 
since for values of a larger than $1$ the solution of Eq.~\eqref{eq:fde-first-order} 
is not monotonic and negative values for $q$ may appear, therefore, it is
not applicable in pharmacokinetics and pharmacodynamics.
Based on these elementary equations the basic relations for the time evolution in 
drug disposition can be formulated, with the assumption of diffusion of drug 
molecules taking place in heterogeneous space. 
The simplest fractional pharmacokinetic model is the one-compartment model with \textit{i.v.} 
bolus administration and the concentration, \textit{c}, can be expressed by
Eq.~\eqref{eq:fde-first-order-solution}  divided by a volume of distribution, as
\begin{align}
\label{eq:fpk-basic}
 c(t) = c_0 \E_{\alpha}(-k_{1,f} t^\alpha),
\end{align}
with $\alpha \leq 1$ and $c_0$ is the ratio [dose]/[apparent volume of distribution]. 
This equation for times $t\ll 1$ behaves as a stretched exponential, i.e., as 
${\sim}\exp(-k_{1,f}t^\alpha/\Gamma(1+\alpha))$, while for large values 
of time it behaves as a power-law, ${\sim}t^{-\alpha}/\Gamma(1-\alpha)$ (see Figure~\ref{fig:f1})~\cite{Mainardi20142267}. 
This kinetics is, therefore, a good candidate to describe various 
datasets exhibiting power-law-like kinetics due to the slow diffusion of the drug in 
deeper tissues.  Moreover, the relevance of the stretched exponential (Weibull) function 
with the Mittag-Leffler function probably explains the successful application of the 
former function in describing drug release in heterogeneous media~\cite{PAPADOPOULOU200644}. 
Eq.~\eqref{eq:fpk-basic} is a relationship for the simplest case of fractional pharmacokinetics. 
It accounts for the anomalous diffusion process, which may be considered to be the 
limiting step of the entire kinetics. Classic clearance may be considered not to be the 
limiting process here and is absent from the equation.

\subsection*{The Laplace transform for FDEs}
Fractional differential equations (FDEs) can be easily written in the Laplace 
domain since each of the fractional derivatives can be transformed similarly to 
the ordinary derivatives, as follows, for order $\alpha \leq 1$:
\begin{align}
 \label{eq:laplace-of-caputo}
 \L\{ {}_{0}^{\mathrm{c}}\D_{t}^{\alpha}f(t)\} = s^\alpha F(s) - s^{\alpha-1}f(0^+),
\end{align}
where $F(s)$ is the Laplace transform of $f(t)$~\cite{Pod99}. 
For $\alpha=1$, Eq.~\eqref{eq:laplace-of-caputo} reduces to the classic well-known
expression for ordinary derivatives, that is, $\L\{\dot{f}(t)\}=sF(s) - f(0^+)$.

Let us take as an example the following simple FDE
\begin{align}
 {}_{0}^{\mathrm{c}}\D_{t}^{1/2}q(t) = -q(t),
\end{align}
with initial value $q(0) = 1$ which can be written in the Laplace domain, 
by virtue of Eq.~\eqref{eq:laplace-of-caputo}, as follows
\begin{align}
 s^{\half}Q(s) - s^{-\half}q(0) = -Q(s),
\end{align}
where $Q(s)=\L\{q(t)\}(s)$ is the Laplace transform of $q$. 
After simple algebraic manipulations, we obtain
\begin{align}
 \label{eq:Qs-simple-example}
 Q(s) = \frac{1}{s+\sqrt{s}}.
\end{align}
By applying the inverse Laplace transform to Eq.~\eqref{eq:Qs-simple-example}
the closed-form analytical solution of this FDE can be obtained; this involves a Mittag-Leffler 
function of order one half. In particular, 
\begin{align}
 q(t) = \E_{\half}(-t^{\half}) = e^t(1+\erf(-\sqrt{t})).
\end{align}

Although more often than not it is easy to transform an FDE in the Laplace domain, 
it is more difficult to apply the inverse Laplace transform so as to 
solve it explicitly for the system variables
so that an analytical solution in the time domain is obtained. However, 
it is possible to perform that step numerically using a numerical inverse 
Laplace transform (NILT) algorithm~\cite{art12} as described above.

\section*{Fractional Pharmacokinetics}

\subsection*{Multi-compartmental models}
A one-compartment pharmacokinetic model with \textit{i.v.} bolus 
administration can be easily fractionalised as 
in Eq.~\eqref{eq:fde-first-order} by changing the derivative 
on the left-hand-side of the single ODE to a fractional order. 
However, in pharmacokinetics and other fields where compartmental models are used, 
two or more ODEs are often necessary and it is not as straightforward 
to fractionalise systems of differential equations, especially when 
certain properties such as mass balance need to be preserved.

When a compartmental model with two or more compartments is 
built, typically an outgoing mass flux becomes an incoming flux to the 
next compartment. Thus, an outgoing mass flux that is defined as a rate 
of fractional order, cannot appear as an incoming flux into another 
compartment, with a different fractional order, without violating mass 
balance~\cite{art14}. It is therefore, in general, impossible to fractionalise 
multi-compartmental systems simply by changing the order of the derivatives on 
the left hand side of the ODEs. 

One approach to the fractionalisation of multi-compartment models 
is to consider a common fractional order that defined the mass transfer 
from a compartment $i$ to a compartment $j$: the outflow of one compartment 
becomes an inflow to the other. This implies a common fractional order known 
as a \textit{commensurate} order.
In the general case, of non-commensurate-order systems, a different approach for fractionalising systems of
ODEs needs to be applied.

In what follows, a general form of a fractional two-compartment system 
is considered and then generalised to a system of an arbitrary number of compartments, 
which first appeared in~\cite{DokMagMac10}. A general ordinary linear two-compartment model, 
is defined by the following system of linear ODEs,
\begin{subequations}
\label{eq:original-int-pk-model}
 \begin{align}
  \frac{\d q_1(t)}{\d t} &= -k_{12}q_1(t) + k_{21}q_2(t) - k_{10}q_1(t) + u_1(t),\\
  \frac{\d q_2(t)}{\d t} &= k_{12}q_1(t) - k_{21}q_2(t) - k_{20}q_2(t) + u_2(t),
 \end{align}
\end{subequations}
where $q_1(t)$ and $q_2(t)$ are the mass or molar amounts of material 
in the respective compartments and the $k_{ij}$ constants control the mass transfer 
between the two compartments and elimination from each of them. 
The notation convention used for the indices of the rate constants is 
that the first corresponds to the source compartment and the second to the target 
one, e.g. $k_{12}$ corresponds to the transfer from compartment 1 to 2, 
$k_{10}$ corresponds to the elimination from compartment 1, etc. 
The units of all the $k_{ij}$ rate constants are (1/[time]). 
$u_{i}(t)$ are input rates in each compartment which may be zero, 
constant or time dependent.  Initial values for $q_1$ and $q_2$ have to be 
considered also, $q_1(0)$ and $q_2(0)$, respectively.

In order to fractionalise this system, first the ordinary system is integrated, 
obtaining a system of integral equations and then the integrals are fractionalised 
as shown in~\cite{DokMagMac10}. Finally, the fractional integral equations, are differentiated 
in an ordinary way. The resulting fractional system contains ordinary derivatives 
on the left hand side and fractional Riemann-Liouville derivatives on the right hand side:
\begin{subequations}
\label{eq:fractionalised-model}
 \begin{align}
  \frac{\d q_1(t)}{\d t} &= 
      - k_{12,f}{\,}_{0}\D_{t}^{1-\alpha_{12}} q_1(t) 
      + k_{21,f}{\,}_{0}\D_{t}^{1-\alpha_{21}} q_2(t) 
      - k_{10,f}{\,}_{0}\D_{t}^{1-\alpha_{10}} q_1(t) + u_1(t),\\
  \frac{\d q_2(t)}{\d t} &= 
        k_{12,f}{\,}_{0}\D_{t}^{1-\alpha_{12}} q_1(t) 
      - k_{21,f}{\,}_{0}\D_{t}^{1-\alpha_{21}}q_2(t) 
      - k_{20,f}{\,}_{0}\D_{t}^{1-\alpha_{20}} q_2(t) + u_2(t),
 \end{align}
\end{subequations}
where $0 < \alpha_{ij} \leq 1$ is a constant representing the order of the specific process. 
Different values for the orders of different processes may be considered, but the order 
of the corresponding terms of a specific process are kept the same when these appear 
in different equations, e.g., there can be an order $\alpha_{12}$ for the transfer 
from compartment 1 to 2 and a different order $\alpha_{21}$ for the transfer from 
compartment 2 to 1, but the order for the corresponding terms of the transfer, from 
compartment 1 to 2, $\alpha_{12}$, is the same in both equations. Also the index $f$ 
in the rate constant was added to emphasise the fact that these are different to the ones 
of Eq.~\eqref{eq:original-int-pk-model} and carry units [time]$^{-\alpha}$.

It is convenient to rewrite the above system Eq.~\eqref{eq:fractionalised-model} with 
Caputo derivatives. An FDE with Caputo derivatives accepts the usual type of 
initial conditions involving the variable itself, as opposed to RL derivatives 
which involve an initial condition on the derivative of the variable, which is 
not practical. When the initial value of $q_1$ or $q_2$ is zero then the respective
RL and  Caputo derivatives are the same. This is convenient since a zero initial 
value is very common in compartmental analysis. When the initial value is not zero, 
converting to a Caputo derivative is possible, for the particular term with a non-zero 
initial value. The conversion from a RL to a Caputo derivative of the form that appears 
in Eq.~\eqref{eq:fractionalised-model} is done with the following expression:
\begin{align}
 \label{eq:rl-to-caputo}
 {}_{0}\D_{t}^{1-\alpha_{ij}} q_i(t)  = 
   {}_{0}^{\mathrm{c}}\D_{t}^{1-\alpha_{ij}}q_i(t) + \frac{q_i(0) t^{\alpha_{ij}-1}}{\Gamma(\alpha_{ij})}.
\end{align}

Summarising the above remarks about initial conditions, we may identify three cases: (i) The initial condition
is zero and then the derivative becomes a Caputo by definition. (ii) The initial condition is 
non-zero, but it is involved in a term with an ordinary derivative so it is treated as usual. 
(iii) The initial condition is non-zero and is involved in a fractional derivative which 
means that in order to present a Caputo derivative, an additional term, involving the initial 
value appears, by substituting Eq.~\eqref{eq:rl-to-caputo}. Alternatively, a zero initial value for that variable
can be assumed, with  a Dirac delta input to account for the initial quantity for that variable. 
So, a general fractional model with two compartments, Eq.~\eqref{eq:fractionalised-model}, was formulated, where the 
fractional derivatives can always be written as Caputo derivatives. It is easy to generalise the
above approach to a system with an arbitrary number of $n$ compartments as follows

\begin{align}
 \label{eq:x234123}
 \frac{\d q_i(t)}{\d t} = 
    - k_{i0} {\,}_{0}^{\mathrm{c}}\D_{t}^{1-\alpha_{i0}}q_i(t)
    - \sum_{j\neq i}k_{ij} {\,}_{0}^{\mathrm{c}}\D_{t}^{1-\alpha_{ij}}q_i(t)
    + \sum_{j\neq i}k_{ji} {\,}_{0}^{\mathrm{c}}\D_{t}^{1-\alpha_{ji}}q_j(t)
    + u_i(t),
\end{align}
for $i=1,\ldots, n$, where Caputo derivatives have been considered throughout since, 
as explained above, this is feasible. This system of Eqs.~\eqref{eq:x234123} is too general 
for most purposes as it allows every compartment to be connected with every other. 
Typically the connection matrix would be much sparser than that, with most compartments 
being connected to just one neighbouring compartment while only a few ``hub'' compartments 
would have more than one connections. 
The advantage of the described approach of fractionalisation is that each transport process 
is fractionalised separately, rather than fractionalising each compartment or each equation. 
Thus, processes of different fractional orders can co-exist since they have consistent 
orders when the corresponding terms appear in different equations. Also, it is important to note 
that dynamical systems of the type~\eqref{eq:x234123} do not suffer from pathologies
such as violation of mass balance or inconsistencies with the units of the rate constants.

As mentioned, FDEs can be easily written in the Laplace domain. In the case of FDEs of the 
form of Eq.~\eqref{eq:x234123}, where the fractional orders are $1-\alpha_{ij}$, the Laplace transform 
of the Caputo derivative becomes
\begin{align}
 \L\{{}_{0}^{\mathrm{c}}\D_{t}^{1-\alpha_{ij}}q_i(t)\} = s^{1-\alpha_{ij}}Q_i(s) - s^{-\alpha_{ij}}q_i(0).
\end{align}
An alternative approach for fractionalisation of non-commensurate fractional pharmacokinetic systems has been proposed in~\cite{POPOVIC20132507}, where the conditions that the pharmacokinetic parameters need to satisfy for the mass balance to be preserved have been defined.

\subsection*{A one-compartment model with constant rate input}
After the simplest fractional pharmacokinetic model with one-compartment and 
\textit{i.v. bolus} of Eq.~\eqref{eq:fpk-basic}, the same model with fractional 
elimination, but with a constant rate input is considered~\cite{DokMagMac10}. 
Even in this simple one-compartment model, it is necessary to 
employ the approach of fractionalising each process separately, described above, since the 
constant rate of infusion is not of fractional order. That would have been difficult if one 
followed the approach of changing the order of the derivative of the left hand side of the ODE, 
however here it is straightforward.
The system can be described by the following equation 
\begin{align}
 \label{eq:single-comp-fde}
 \frac{\d q(t)}{\d t} = k_{01} - k_{10,f}{\,}_{0}^{\mathrm{c}}\D_{t}^{1-\alpha}q(t),
\end{align}
with $q(0)=0$ and where $k_{01}$ is a zero-order input rate constant, 
with units [mass]/[time], $k_{10,f}$ is a rate constant with units 
[time]$^{-\alpha}$ and $\alpha$ is a fractional order less than $1$.
Eq.~\eqref{eq:single-comp-fde} can be written in the Laplace domain as
\begin{align}
\label{eq:single-comp-fde-s}
 sQ(s) - q(0) = \frac{k_{01}}{s} - k_{10,f}(s^{1-\alpha}Q(s) - s^{-\alpha}q(0)).
\end{align}
Since $q(0)=0$, Eq.~\eqref{eq:single-comp-fde-s} can be solved to obtain
\begin{align}
 Q(s) = \frac{k_{01}s^{\alpha-2}}{s^{\alpha}+k_{10,f}}.
\end{align}
By applying the following inverse Laplace transform formula (Equation (1.80) in~\cite{Pod99}, page 21):
\begin{align}
 \label{eq:laplace-two-param-MLf}
 \L^{-1}\left\{ \frac{s^{\mu-\nu}}{s^{\mu} + k}\right\} = t^{\nu-1}\E_{\mu,\nu}(-k t^\mu),
\end{align}
where $\E_{\mu,\nu}$ is the Mittag-Leffler function with two parameters. For $\mu=\alpha$ and $\nu = 2$,
we obtain
\begin{align}
 \label{eq:single-comp-fde-sol}
 q(t) = k_{01} t \E_{\alpha,2}(-k_{10,f}t^{\alpha}).
\end{align}

In Theorem 1.4 of~\cite{Pod99}, the following expansion for the Mittag-Leffler 
function is proven to hold asymptotically as $|z|\to\infty$:
\begin{align}
 \E_{\mu,\nu}(z) = -\sum_{k=1}^{p}\frac{z^{-k}}{\Gamma(\nu-\mu)} + \mathcal{O}(|z|^{1-p}).
\end{align}
Applying this formula for Eq.~\eqref{eq:single-comp-fde-sol} and keeping only 
the first term of the sum since the rest are of higher order, the limit 
of Eq.~\eqref{eq:single-comp-fde-sol} for $t$ going to infinity can be shown that it is $\lim_{t\to\infty}q(t)=\infty$ for all $0 < \alpha < 1$~\cite{DokMagMac10}.

\begin{figure}
 \centering
 \includegraphics[width=0.6\textwidth]{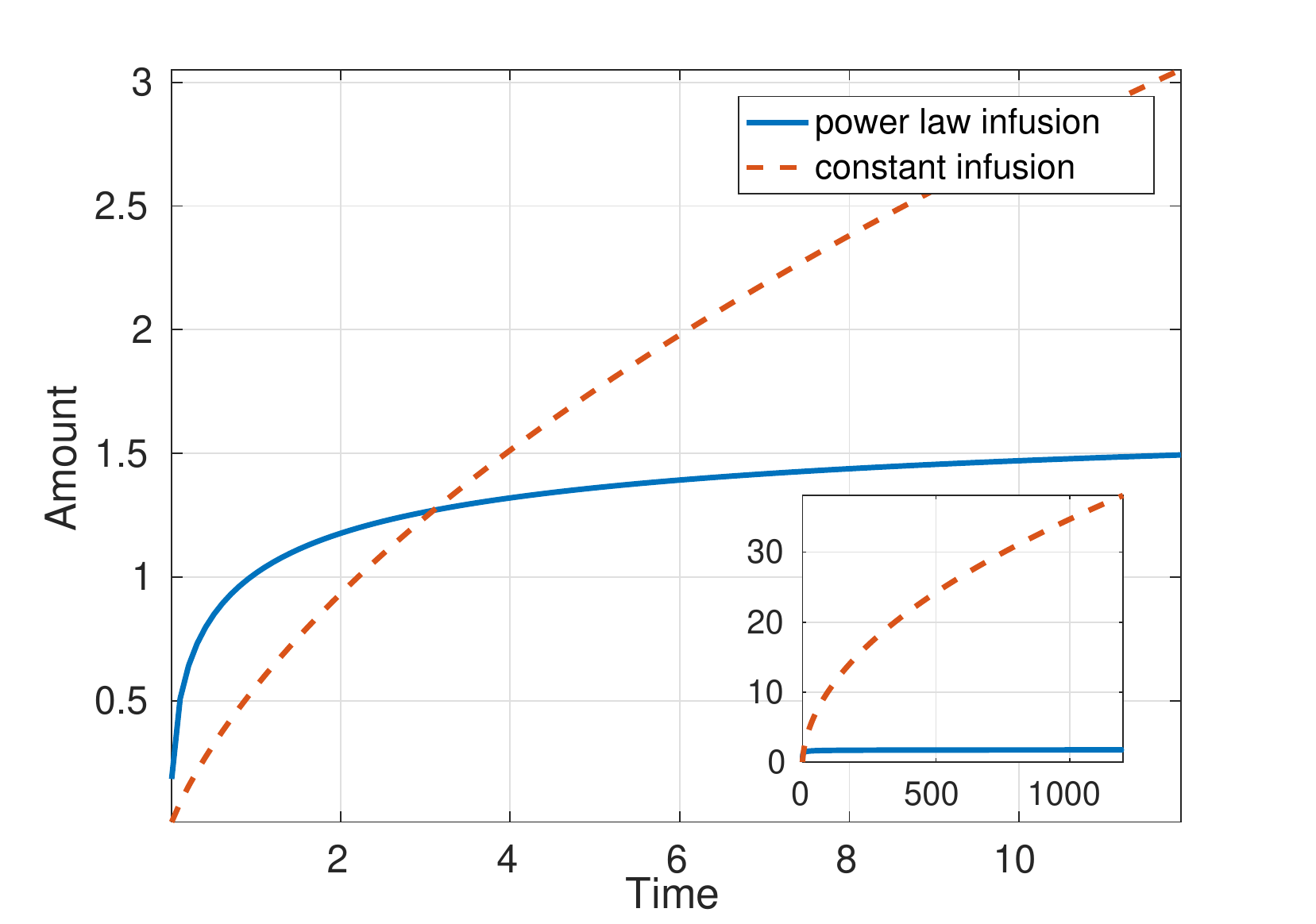}
 \caption{Amount-time profiles for $\alpha=0.5$: (solid line)
          amount of drug versus time according to Eq.~\eqref{eq:single-comp-fde-sol} 
          with constant infusion where there is unbounded accumulation of drug, 
          (dashed line) time-profile according to Eq.~\eqref{eq:single-comp-fde-sol2}, 
          with power-law infusion where the amount of drug approaches a steady state. 
	  (Inset) The same profiles for $100$ times longer time span.}
 \label{fig:f2}	  
\end{figure}

The fact that the limit of $q(t)$ in Eq.~\eqref{eq:single-comp-fde-sol} 
diverges as $t$ goes to infinity, 
for $\alpha<1$, means that unlike the classic case, for $\alpha=1$ --- where 
\eqref{eq:single-comp-fde-sol} approaches exponentially the steady state 
$\nicefrac{k_{01}}{k_{10,f}}$, for $\alpha<1$ --- there is infinite accumulation. 
In Figure~\ref{fig:f2} a plot of~\eqref{eq:single-comp-fde-sol} is shown for $\alpha=0.5$ 
demonstrating that in the fractional case the amount grows unbounded. 
In the inset of Figure~\ref{fig:f2} the same profiles are shown a 100-fold larger time span, 
demonstrating the effect of continuous accumulation.

The lack of a steady state under constant rate administration which results 
to infinite accumulation is one of the most important clinical implications 
of fractional pharmacokinetics. It is clear that this 
implication extents to repeated doses as well as constant infusion, which is the 
most common dosing regimen, and can be important in chronic treatments. 
In order to avoid accumulation, the constant rate administration must be adjusted 
to a rate which decreases with time. Indeed, in Eq.~\eqref{eq:single-comp-fde}, if the constant rate of 
infusion $k_{01}$ is replaced by the term $f(t)=k_{01}t^{1-\alpha}$,~\cite{Hennion2013} 
the solution of the resulting FDE is, instead of Eq.~\eqref{eq:single-comp-fde-sol}, 
the following
\begin{align}
 \label{eq:single-comp-fde-sol2}
 q(t) = k_{01}\Gamma(\alpha) t^\alpha \E_{\alpha, \alpha+1}(-k_{10,f}t^{\alpha}).
\end{align}

The drug quantity $q(t)$ in Eq.~\eqref{eq:single-comp-fde-sol2} converges to the steady 
state $\Gamma(\alpha)\nicefrac{k_{01}}{k_{10,f}}$ 
as time goes to infinity, while for the special case of $\alpha=1$, the steady 
state is the usual $\nicefrac{k_{01}}{k_{10,f}}$. Similarly, for the case of repeated doses, 
if a steady state is intended to be achieved, in the presence of fractional elimination 
of order $\alpha$, then the usual constant rate of administration, e.g., a constant daily dose, 
needs to be replaced by an appropriately decreasing rate of administration. As shown in~\cite{Hennion2013}, 
this can be either the same dose, but given at increasing inter-dose intervals, i.e.,
$T_i = (T_{i-1}^{\alpha} + \alpha \Delta \tau^{\alpha})^{\nicefrac{1}{\alpha}}$, where
$T_i$ is the time of the $i$-th dose and $\Delta \tau$ is the inter-dose interval of the 
corresponding kinetics of order $\alpha=1$; 
or a decreasing dose given at constant intervals,
i.e., $q_{0,i}=\nicefrac{q_0}{\alpha}((i+1)^{\alpha}-i^{\alpha})$. In this way, an ever 
decreasing administration rate is implemented which compensates the decreasing elimination 
rate due to the fractional kinetics. 

\subsection*{A two-compartment i.v. model}
\begin{figure}
 \centering
 \includegraphics[width=0.55\textwidth]{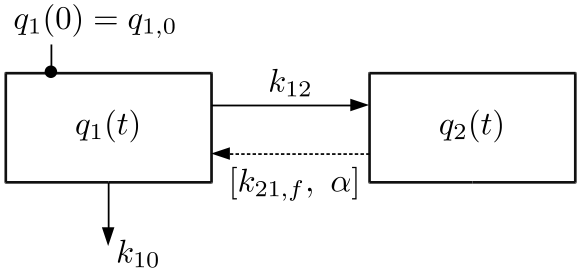}
 \caption{A fractional two-compartment PK model with an \textit{i.v.} bolus. Elimination from the 
	  central compartment and a mass flux from the central to the peripheral compartment, 
	  which are both assumed to follow classic kinetics (order 1), while a flux from the 
	  peripheral to the central compartment is assumed to follow slower fractional 
	  kinetics, accounting for tissue trapping (dashed arrow).}
 \label{fig:f3}	  
\end{figure}

Based on the generalised approach for the fractionalisation of compartmental models, 
which allows mixing different fractional orders, developed above, 
a two-compartment fractional pharmacokinetic model is considered, shown schematically in 
Figure~\ref{fig:f3}. Compartment 1 (central) represents general circulation and well perfused tissues 
while compartment 2  (peripheral) represents deeper tissues. Three transfer processes (fluxes) 
are considered: elimination from the central compartment and a mass flux from the central to 
the peripheral compartment, which are both assumed to follow classic kinetics (order 1), while 
a flux from the peripheral to the central compartment is assumed to follow slower fractional 
kinetics accounting for tissue trapping~\cite{Petras2011,art13}.

The system is formulated mathematically as follows:
\begin{subequations}
\label{eq:model-two-comp}
 \begin{align}
  \frac{\d q_1(t)}{\d t} &= -(k_{12} + k_{10})q_1(t) 
     + k_{21,f} {\,}_{0}^{\mathrm{c}}\D_{t}^{1-\alpha}q_2(t),\\
  \frac{\d q_2(t)}{\d t} &=  k_{12}q_1(t) - k_{21,f}
          {\,}_{0}^{\mathrm{c}}\D_{t}^{1-\alpha}q_2(t),
 \end{align}
\end{subequations}
where $\alpha < 1$ and initial conditions are $q_1(0) = q_{1,0}$, $q_2(0) = 0$ 
which account for a bolus dose injection and no initial amount in the peripheral 
compartment. Note, that it is allowed to use Caputo derivatives here since the 
fractional derivatives involve only terms with $q_2$ for which there is no 
initial amount, which means that Caputo and RL derivatives are identical.

Applying the Laplace transform, to the above system the following algebraic equations are obtained:
\begin{subequations}
 \begin{align}
  sQ_1(s) - q_1(0) &= - (k_{12} + k_{10})Q_1(s) + k_{21}(s^{1-\alpha}Q_2(s) - s^{-\alpha}q_2(0)),\\
  sQ_2(s) - q_2(0) &= k_{12}Q_1(s) - k_{21,f}(s^{1-\alpha}Q_2(s) - s^{-\alpha}q_2(0)).
 \end{align}
\end{subequations}

Solving for $Q_1(s)$ and $Q_2(s)$ and substituting the initial conditions, 
\begin{subequations}
 \begin{align}
  Q_1(s) &= 
      \frac
	  {q_{1,0}(s^{\alpha} + k_{21,f})}
	  {(s+k_{12}+k_{10})(s^\alpha+k_{21,f})-k_{12}k_{21,f}}
	  \label{eq:Q1s-two-compartment}\\
  Q_2(s) &= 
      \frac
	  {q_{1,0}  s^{\alpha-1} k_{12}} 
	  {(s+k_{12}+k_{10})(s^\alpha+k_{21,f})-k_{12}k_{21,f}}	  
	  \label{eq:Q2s-two-compartment}
 \end{align}
\end{subequations}

Using the above expression for $Q_1$ and $Q_2$, Eq.~\eqref{eq:Q1s-two-compartment} and Eq.~\eqref{eq:Q2s-two-compartment}, 
respectively, can be used to simulate values of $q_1(t)$ and $q_2(t)$ 
in the time domain by an NILT method. Note, that primarily, $q_1(t)$ is 
of interest, since in practice we only have data from this compartment. 
The output for $q_1(t)$ from this numerical solution may be combined with the 
following equation 
\begin{align}
 \label{eq:conc-quantity}
 c(t) = \frac{q_1(t)}{V},
\end{align}
where $c$ is the drug concentration in the blood and $V$ is the apparent volume of 
distribution. Eq.~\eqref{eq:conc-quantity} can be fitted to pharmacokinetic data in 
order to estimate parameters $V$, $k_{10}$, $k_{12}$, $k_{21,f}$ and $\alpha$.

The closed-form analytical solution of Eq.~\eqref{eq:model-two-comp}, can be expressed in terms 
of an infinite series of generalised Wright functions as demonstrated in the 
book by Kilbas \textit{et al.}~\cite{book17}, but these solutions are hard to implement and apply 
in practice. Moloni in~\cite{art18} derived analytically the inverse Laplace transform of $Q_1(s)$ as:
\begin{subequations}
\label{eq:complex-solution}
 \begin{align}
  q_1(t) &= q_{1,0} 
      \sum_{n=0}^{\infty}(-1)^n k_{21,f}^n 
      \sum_{l=0}^{n} \frac
		  {k_{10}^l n!}
		  {(n-l)!l!}      
      \Big( t^{l+\alpha n}\E_{1,l+\alpha n + 1}^{n+1}
         (-(k_{10}+k_{12})t)  \notag\\
         &+t^{l+\alpha (n+1)} \E_{1,l+\alpha (n+1) + 1}^{n+1}
         (-(k_{10}+k_{12})t)  \Big),
         \label{eq:complex-solution-1}
 \end{align}
while for $Q_2$ the inverse Laplace transform works out to be~\cite{art18}:
\begin{align}
 q_2(t) = q_{1,0} k_{12}
   \sum_{n=0}^{\infty}(-1)^n k_{21,f}^n 
   \sum_{l=0}^{n} \frac
		  {k_{10}^l n!}
		  {(n-l)!l!}
  t^{l+\alpha n}
  \E_{1,l+\alpha n + 2}^{n+1}
   (-(k_{10}+k_{12})t).
\end{align}
\end{subequations}

\begin{figure}[ht!]
 \centering
 \includegraphics[width=0.6\textwidth]{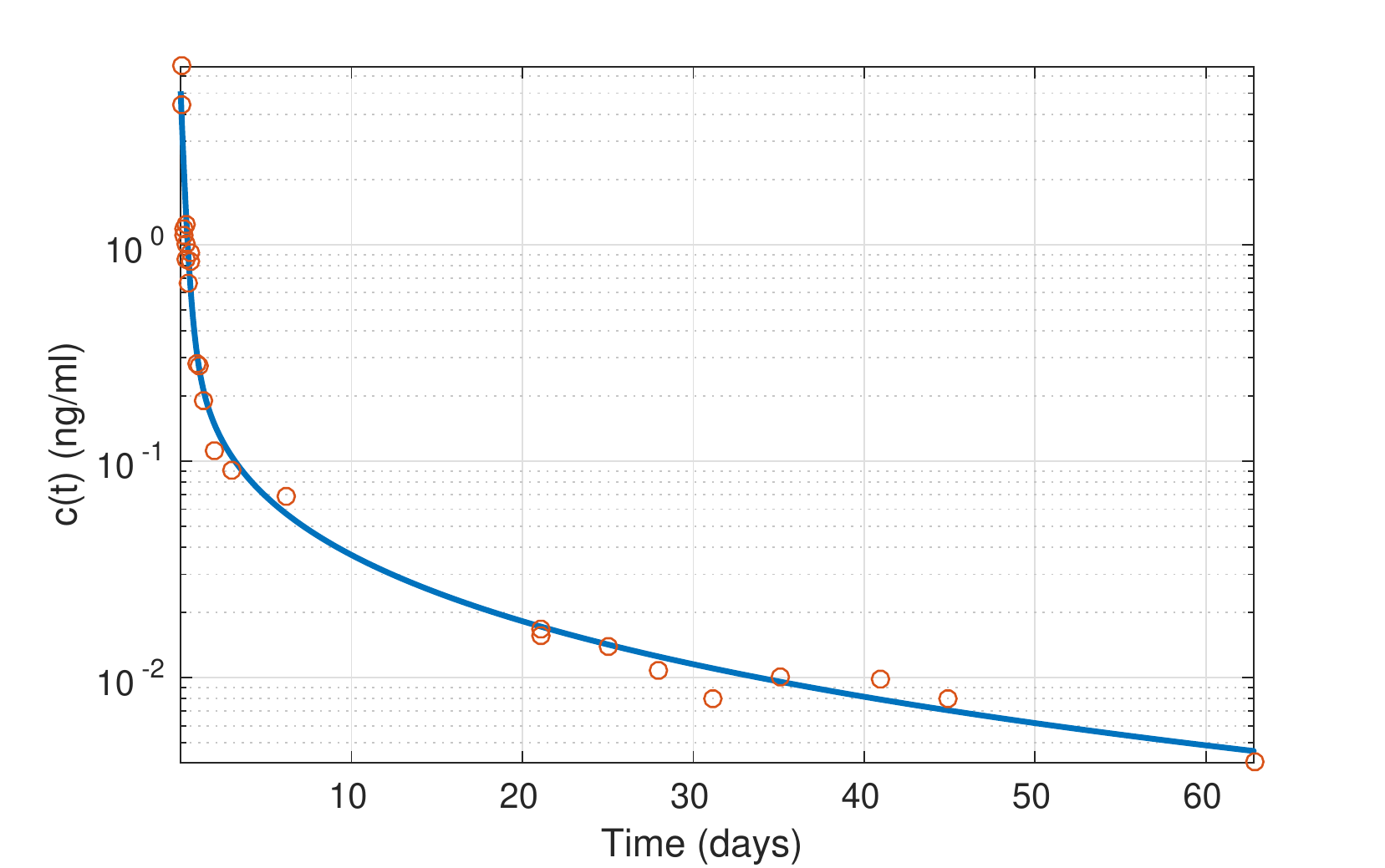}
 \caption{Concentration-time profile of amiodarone in the central compartment.
          The solid line corresponds to the fractional-order pharmacokinetic 
          model of Eq.~\eqref{eq:model-two-comp} with the parameter estimates: 
          $\alpha=0.587$,
	  $k_{10}=\unit[1.49]{days}^{-1}$, 
	  $k_{12}=\unit[2.95]{days}^{-1}$, 
	  $k_{21,f}=\unit[0.48]{days}^{-\alpha}$ and
	  $\nicefrac{q_0}{V}=\unitfrac[4.72]{ng}{ml}$.
	  	  The circles correspond to experimental measurements.}%
 \label{fig:f4}
\end{figure}

An application of the fractional two compartment model,
the system of Eq.~\eqref{eq:model-two-comp} to 
amiodarone PK data was presented in~\cite{DokMagMac10}. 
Amidarone is an antiarrhythmic drug known for 
its anomalous, non-exponential pharmacokinetics, which has important clinical 
implications due to the accumulation pattern of the drug in long-term administration. 
The fractional two-compartment model of the previous section was used to analyse an 
amiodarone \textit{i.v.} data-set which first appeared in~\cite{art19} and estimates of the model 
parameters were obtained. Analysis was carried out in  MATLAB while the values 
for $q_1(t)$ of Eq.~\eqref{eq:complex-solution-1} were simulated using a NILT algorithm~\cite{art12}
from the expression of $Q_1(s)$ in the Laplace domain. 
In Figure~\ref{fig:f4} the model-based predictions are plotted together the data demonstrating 
good agreement for the 60 day period of this study. 
The estimated fractional order was $\alpha=0.587$ and the non-exponential character of 
the curve is evident, while the model follows well the data both for long and 
for short times, unlike empirical power-laws which explode at $t=0$.

\section*{Numerical methods for fractional-order systems}%
For simple fractional-order models, as we discussed previously,
there may exist closed-form analytical solutions which involve the one-parameter
or two-parameters Mittage-Leffler function, or are given by more 
intricate analytical expressions such as Eq.~\eqref{eq:complex-solution}.
Interestingly, even for simple analytical solutions such as 
Eq.~\eqref{eq:fde-first-order-solution}, the Mittag-Leffler function
itself is evaluated by solving an FDE numerically~\cite{Kaczorek2011}.
This fact, without discrediting the value of analytical solutions,
necessitates the availability of reliable numerical methods 
that allow us to simulate and study fractional-order systems.
The availability of accurate discrete-time approximations
of the trajectories of such systems is important not only 
for simulating, but also for the design of 
open-loop or closed-loop administration strategies, based on control theory~\cite{SopNtoSar15}.
Time-domain approximations are less parsimonious than $s$-domain ones,
but are more suitable for control applications as we discuss in 
the next section.

There can be identified four types of solutions for
fractional-order differential equations: (i) closed-form analytical
solutions, (ii) approximations in the Laplace $s$-domain using 
integer-order rational transfer functions, (iii) numerical 
approximation schemes in the discrete time domain and (iv)
numerical inverse Laplace transforms.
Each of these comes with certain advantages as well as limitations; for example,
closed-form analytical solutions are often not available, while the inverse 
Laplace function requires an explicit closed form for $u_i(t)$
so it cannot be used for administration rates that are defined implicitly
or are arbitrary signals. Approximations in the $s$-domain are powerful
modelling tools, but they fail to provide error bounds on the concentration 
of the administered drug in the time domain which are necessary 
in clinical practice. 

In regard to closed-form analytical solutions, when available,
they involve special functions such as the Mittag-Leffler 
function $\mathcal{E}_{\alpha,\beta}(t)=\sum_{k=0}^{\infty}
{t^k}/{\Gamma(\alpha k + \beta)}$ whose evaluation 
calls, in turn, for some numerical approximation scheme. 
Analytical closed-form solutions of fractional differential 
equations are available for pharmacokinetic systems~\cite{Verotta2010}.
Typically for the evaluation of this function we resort 
to solving an FDE numerically~\cite{Garrappa15,Seybold2009,Gorenflo2002}
because the series in the definition of $\mathcal{E}_{\alpha,\beta}(t)$ converges
rather slowly and no exact error bounds as available so as to establish 
meaningful termination criteria.

\subsection*{Transfer functions and integer-order approximations}
Fractional-order systems, like their integer-order counterparts,
can be modelled in the Laplace $s$-domain via \textit{transfer functions},
that is, 
\begin{align}
 F_{ij}(s) = \frac{Q_j(s)}{U_i(s)},
\end{align}
where $Q_j(s)$ and $U_i(s)$ are the Laplace transforms of the drug quantity
$q_j(t)$ and the administration rate $u_i(t)$. If the pharmacokinetics
is described by linear fractional-order models such as the ones discussed 
above, $F_{ij}$ will be a fractional-order rational function (a quotient
of polynomials with real exponents).

Rational approximations aim at approximating such transfer
functions --- which involve terms of the form $s^\alpha$, $\alpha\in\Re$ 
--- by ordinary transfer functions of the form 
\begin{align}
 \tilde{F}_{ij}(s) = \frac{P_{ij}(s)}{S_{ij}(s)},
\end{align}
where $P_{ij}$ and $S_{ij}$ are polynomials and the degree of $P_{ij}$ is
no larger than the degree of $S_{ij}$.
For convenience with notation, we henceforth drop the subscripts $ij$.

%
% ---- PADE ----
%
\textit{Pad\'{e} Approximation:} 
The Pad\'{e} approximation of order $[m/n]$, $m,n\in\N$,
at a point $s_0$ is rather popular and leads to rational functions 
with $\deg P =m$ and $\deg S = n$~\cite{SILVA2006373}. 

%
% ---- MATSUDA-FUJII ----
%
\textit{Matsuda-Fujii continuous fractions approximation:}
This method consists in interpolating a function $F(s)$, which is treated
as a black box, across a set of logarithmically spaced points~\cite{Matsuda1993}. 
By letting the selected points be $s_k$, $k=0,1,2,\dots$, the approximation is 
written as the continued fractions expansion
\begin{equation}
F(s) = \alpha_0 + \frac{s-s_0}{\alpha_1 + \frac{s-s_1}{\alpha_2 + \frac{s-s_2}{\alpha_3+\dots}}}
\end{equation}
where, $\alpha_i = \upsilon_i (s_i)$,  $\upsilon_0(s) = H(s)$ and $\upsilon_{i+1}(s) 
= \frac{s-s_i}{\upsilon_i(s) - \alpha_i}$.

%
% ---- OUSTALOUP ----
%
\textit{Oustaloup's method:}
Oustaloup's method is based on the approximation of a function of the form:
\begin{equation}
H(s) = s^{\alpha}, 
\end{equation}
with $\alpha>0$ by a rational function
\begin{equation}\label{eq:oustaloup}
\widehat{H}(s) = c_0 \prod_{k=-N}^N \frac{s+\omega_k}{s+\omega'_k},
\end{equation}
within a range of frequencies from $\omega_b$ to $\omega_h$~\cite{OusLev+00}.
The Oustaloup method offers an approximation at frequencies 
which are geometrically distributed about the characteristic 
frequency $\omega_u = \sqrt{\omega_b\omega_h}$ --- the geometric 
mean of $\omega_b$ and $\omega_h$.
The parameters $\omega_k$ and $\omega_k'$ are determined via the design 
formulas~\cite{Petras11}
\begin{subequations}
 \begin{align}
  \omega_k' = \omega_b \left(\frac{\omega_h}{\omega_b}\right)^{
		\frac{k+N+0.5(1+\alpha)}{2N+1}
	      },\\
  \omega_k = \omega_b \left(\frac{\omega_h}{\omega_b}\right)^{
		\frac{k+N+0.5(1-\alpha)}{2N+1}
	      },\\
  c_0 = \left(\frac{\omega_h}{\omega_b}\right)^{-\frac{r}{2}}
        \prod_{k=-N}^{N}\frac{\omega_k}{\omega_k'}.
 \end{align}
\end{subequations}
Parameters $\omega_b$, $\omega_h$ and $N$ are design parameters 
of the Oustaloup method.

%
% ---- OTHER METHODS ----
%  ~ Charef et al.
%  ~ Carlson-Halijak
%  ~ System ID methods
%
There are a few more methods which have been proposed in the 
literature to approximate fractional-order systems by 
rational transfer functions such as~\cite{Charef,CarHal64}, 
as well as data-driven system identification techniques~\cite{Gao2012}.
%Nice~\cite{state-space-frac}.

\subsection*{Time domain approximations}
Several methods have been proposed which attempt to 
approximate the solution to a fractional-order initial
value problem in the time domain.

%
% ---- GRUNWALD-LETNIKOV ----
%
\textit{Gr\"{u}nwald-Letnikov:}
This is the method of choice in the discrete time
domain where $\glD^\alpha x(t)$ is approximated 
by its discrete-time finite-history variant
$
 ({\glDD}^\alpha_{h,\nu} x)_k
$
which is proven to have bounded error with respect to 
$({\glDD}_h^\alpha x)_k$~\cite{SopSar16}.
The boundedness of this approximation error is a singular
characteristic of this approximation method and is suitable
for applications in drug administration where it is necessary
to guarantee that the drug concentration in the body will
remain within certain limits~\cite{SopSar16}.
As an example, system~\eqref{eq:fde-first-order} is approximated (with sampling time $h$) 
by the discrete-time linear system
\begin{align}
 \tfrac{1}{h^{\alpha}}
   \sum_{j=0}^{\min\{\nu,\left\lfloor \nicefrac{t}{h}\right\rfloor\}}
   c_j^\alpha q_{k+1-j} = -k_{1,f} q_{k},
\end{align}
where $q_k = q(kh)$.

The discretisation of the Gr\"{u}nwald-Letnikov derivative suffers from the fact that 
the required memory grows unbounded with the simulation time. The truncation of the Gr\"{u}nwald-Letnikov 
series up to some finite memory gives rise to viable solution alorithms which can 
be elegantly described using \textit{triangular strip matrices} as described in~\cite{Podlubny00matrixapproach}
and are available as a MATLAB toolbox.

\textit{Approximation by parametrisation:}
Approximate time-domain solutions can be obtained by assuming a particular 
parametric form for the solution. 
Such a method was proposed by Hennion and Hanert~\cite{Hennion2013}
where $q(t)$ is approximated by finite-length expansions of the form $\sum_{j=0}^{N}A_j \phi_j(t)$,
where $\phi_j$ are Chebyshev polynomials and $A_j$ are constant coefficients. 
By virtue of the computability of fractional derivatives of $\phi_j$, the parametric 
approximation is plugged into the fractional differential equation
which, along with the initial conditions, yields a linear system.
What is notable in this method is that expansions as short as $N=10$ lead to very 
low approximation errors.
Likewise, other parametric forms can be used. For example,~\cite{Zainal2014}
used a Fourier-type expansion and~\cite{KUMAR20062602} used piecewise 
quadratic functions.

%
% ---- NUMERICAL METHODS in general ----
%
\textit{Numerical integration methods:}
Fractional-order initial value problems can be solved
with various numerical methods such as the Adams-Bashforth-Moulton 
predictor-corrector (ABMPC) method~\cite{Zayernouri20161} and 
fractional linear multi-step methods (FLMMs)~\cite{Lubich,Garrappa_trapezoidal}.
These methods are only suitable for systems of FDEs in the form
\begin{subequations}\label{eq:FDE_basic_form}
\begin{align}
 {}_{0}^{\mathrm{c}}\D_{t}^\gamma x(t) &= f(t,x(t)),\\
 x^{(k)}(0) &= x_{0,k},\, k = 0,\dots,m-1
\end{align}
\end{subequations}
where $\gamma$ is a rational, and $m = \lceil \gamma \rceil$.
Typically in pharmacokinetics we encounter cases with $0 < \gamma \leq 1$,
therefore, we will have $m=1$.

Let us give an example on how this applies to the two-compartment model we 
presented above. 
In order to bring Eq.~\eqref{eq:model-two-comp} in the above form, 
we need to find a rational approximation of two derivatives, $1-\alpha$ and $1$. 
If we can find a satisfying rational approximation of $1 - a\approxeq p/q$,
then the first order derivative follows trivially.
Now, Eq.~\eqref{eq:model-two-comp} can be written as 
\begin{subequations}\label{eq:pkModelAmiodarone_system}
\begin{align}
    {}_{0}^{\mathrm{c}}\D_{t}^{\gamma} x_0 &= x_1 \\ 
    {}_{0}^{\mathrm{c}}\D_{t}^{\gamma} x_1 &= x_2 \\
    &\ \,\vdots \nonumber \\
    {}_{0}^{\mathrm{c}}\D_{t}^{\gamma} x_{q-1} &= -(k_{12}+k_{10})x_0 {+} k_{21,f} x_{q+p} {+} u\\
    {}_{0}^{\mathrm{c}}\D_{t}^{\gamma} x_{q} &=  x_{q+1}  \\
    &\ \, \vdots \nonumber \\
    {}_{0}^{\mathrm{c}}\D_{t}^{\gamma} x_{2q-1} &=  k_{12}x_0 - k_{21,f}x_{q+p}
\end{align}
\end{subequations}
with initial conditions $x_0(0) = q_1(0)$, $x_{q}(0) = q_2(0)$ and $x_i(0) = 0$ for $i \notin \{1, q\}$,
and $\gamma=1/q$.
This system is in fact a linear fractional-order 
system for which closed-form analytical solutions are available~\cite[Thm.~2.5]{Kaczorek2011}.
In particular, Eq.~\eqref{eq:pkModelAmiodarone_system} can be written concisely in the 
form 
\begin{subequations}
 \begin{align}
  {}_{0}^{\mathrm{c}}\D_{t}^{\gamma} \mathbf{x}(t) &= A\mathbf{x}(t) + Bu(t)\\
  \mathbf{x}(0) &= \mathbf{x}_0,
 \end{align}
\end{subequations}
where $\mathbf{x} = (x_0, x_1,\ldots, x_{2q-1})$ and $\mathbf{x}_0 = (x_0(0),\ldots, x_{2q-1}(0))$
and matrices $A$ and $B$ can be readily constructed from the above dynamical equations.
This fractional-order initial value problem has the closed-form analytical 
solution~\cite[Thm.~2.5]{Kaczorek2011}
\begin{align}
 \mathbf{x}(t) = \E_\gamma (At^\gamma)\mathbf{x}_0 + \int_0^t 
   \sum_{k=0}^{\infty}\frac{A^k(t-\tau)^{k\gamma}}{\Gamma((k+1)\gamma)}Bu(\tau)\d\tau.
\end{align}
Evidently, the inherent complexity of this closed-form analytical solution --- which would
require the evaluation of slowly-converging series --- motivates and necessitates 
the use of numerical methods.

The number of states of system~\eqref{eq:pkModelAmiodarone_system} is $2q$, 
therefore, the rational approximation should 
aim at a small $q$. Yet another reason to choose small $q$ is that small values of 
$\gamma$ %=1/q$% 
render the system hard to simulate numerically because they increase the effect of and 
dependence on the memory.
 %In our case $1-\alpha= 0.413$ can be written 
%as $413/1000$, but then we would need to simulate a fractional-order
%system with $2000$ states.
% A reasonable approximation of $1-\alpha$ is 
% $19/46$ with error $0.413-19/46 = -4.3478\cdot 10^{-5}$. % or 
% %$216/523$ with error $0.413-216/523 = -1.912\cdot 10^{-6}$.
% Such approximations can be obtained by means of continued fractions 
% expansions of $1-\alpha$.

%
% --- ABMPC methods ---
%
\textit{Adams-Bashforth-Moulton predictor-corrector (ABMPC):} 
Methods of the ABMPC type  have been generalised 
to solve fractional-order systems.
The key idea is to evaluate $({\rlI}^\gamma f)(t,x(t))$ by 
approximating $f$ with appropriately selected polynomials.
Solutions of Eq.~\eqref{eq:FDE_basic_form} satisfy the following
integral representation
\begin{align}
 x(t) = \sum_{k = 0}^{m-1} x_{0,k}\frac{t^k}{k!} + ({\rlI}^\gamma f)(t,x(t)),
\end{align}
where the first term on right hand side will be denoted by $T_{m-1}(t).$
This is precisely an integral representation of Eq.~\eqref{eq:FDE_basic_form}.
The integral on the right hand side of the previous equation can 
be approximated, using an uniformly spaced grid $t_n = nh$, by
$
  \frac{h^\gamma}{\gamma(\gamma + 1)}\sum_{j = 0}^{n+1} a_{j,n+1} f(t_j) 
$ for suitable coefficients $a_{j,n+1}$~\cite{Diethelm2002}.
The numerical approximation of the solution of Eq.~\eqref{eq:FDE_basic_form} is 
\begin{subequations}
\begin{align}\label{eq:corrector_equation}
 x(t_{n+1}) &= T_{m-1}(t_{n+1}) 
             + \tfrac{h^\gamma}
                     {\Gamma(\gamma+2)}f(t_{n+1},x_p(t_{n+1}))
            + \sum_{j=1}^n a_{j,n+1}f(t_j,x(t_j)). 
 \end{align}
The equation above is usually referred to as the corrector formula and 
$x_p(t_{n+1})$ is given by the predictor formula
\begin{align}
 x_P(t_{n+1}) = T_{m-1}(t_n) &+ \tfrac{1}{\Gamma(\gamma)} 
                 \sum_{j = 0}^{n} b_{j,n+1}f(t_j,x(t_j)).
\label{eq:predictor_equation}
\end{align}
\end{subequations} 
Unfortunately, the convergence error of ABMPC when $0 < \gamma < 1$
is $\mathcal{O}(h^{1+\gamma})$, therefore, a rather small step size 
$h$ is required to attain a reasonable approximation error.
A modification of the basic predictor-corrector method with more 
favourable computational cost is provided in~\cite{Garrappa_pc_stability} 
for which the MATLAB implementation \texttt{fde12} is available.

\textit{Lubich's method:} Fractional linear multi-step 
methods (FLMM)~\cite{Lubich} are a generalisation of 
linear multi-step methods (LMM) for ordinary differential equations.
The idea is to approximate the 
Riemann-Liouville fractional-order integral operator~\eqref{eq:riemann-liouville-integral}
with a discrete convolution, known as a \textit{convolution quadrature}, as
\begin{align}
 (\rlI_h^\gamma f)(t) 
   \approxeq h^\gamma\sum_{j=0}^{n} \omega_{n-j}f(t_j) 
   + h^\gamma \sum_{j=0}^{s}w_{n,j}f(t_j),
 \label{eq:flmm_conv_quadrature}
\end{align}
for $t_j = jh$, where ($w_{n,j}$) and ($\omega_n$) are independent of $h$. 
%
% (!) Explain...
Surprisingly, the latter can be constructed from any linear 
multistep method for arbitrary fractional order $\gamma$~\cite{Lubich}.
FLMM constructed this way will inherit the same 
convergence rate and at least the same stability properties as the 
original LMM method~\cite{Lubich_flmm_abel_voltera}.

A MATLAB implementation of Lubich's method, namely \texttt{flmm2}~\cite{Garrappa_software} 
which is based on~\cite{Garrappa_trapezoidal} is available.
However, these methods do not perform well for small $\gamma$.
According to~\cite{HerNto+17}, for the case of amiodarone,
it is reported that values of  $\gamma$ smaller than  
$0.1$ give poor results and often do not converge, while using
a crude approximation given by $\gamma = 1/5$, 
\texttt{flmm2} was shown to outperform \texttt{fde12} in terms of accuracy and 
stability at bigger step sizes $h$.

% (!) Comment: gamma is desirable to be not too small!
\subsection*{Numerical inverse Laplace}
The inverse Laplace transform of a transfer function $F(s)$ 
---  on the Laplace $s$-domain --- is given by the complex integral
\begin{align}
  f(t) = \L^{-1}\{F\}(t) = \lim_{T \rightarrow \infty}\frac{1}{2\pi i} \int_{\sigma - iT}^{\sigma + iT} e^{st}F(s)\d s,
  \label{eq:inv_laplace}
\end{align}
where $\sigma$ is any real number larger than the real parts of the 
poles of $F$.
The numerical inverse Laplace (NILT) approach aims at approximating the 
above integral numerically. 
Such numerical methods apply also to cases where $F$ is not rational
as it is the case with fractional-order systems~\eqref{eq:inv_laplace}.

In the approach of De Hoog \textit{et al.}, the integral which 
appears in Eq.~\eqref{eq:inv_laplace} is cast as a Fourier transform
which can then be approximated by a Fourier series followed by 
standard numerical integration (e.g., using the trapezoid rule)~\cite{art12}. 
Though quite accurate for a broad class of functions, these methods are typically
very demanding from the computational point of view.
An implementation of the above method is available online~\cite{Hollenbeck}.

In special cases analytical inversion can be done by means of 
Mittag-Leffler function~\cite{Lin2013,Kexue}.
A somewhat different approach is taken by \cite{valsa_invlap},
where authors approximate $e^{st}$, the kernel of the 
inverse Laplace transformation, by a function of the form $\tfrac{e^{st}}{1 + e^{-2a}e^{2st}}$
and choose $a$ appropriately so as to achieve an accurate inversion.

In general, numerical inversion methods can achieve high precision, but 
they are not suitable for control design purposes, especially
for optimal control problems. 
Moreover, there is not one single method
that gives the most accurate inversion for all types of functions.
An overview of the most popular inversion methods used in 
engineering practice is given in~\cite{invlap_comp}.

\section*{Drug administration for fractional pharmacokinetics}
Approaches for drug administration scheduling can be classified according to 
the desired objective into methods where (i) we aim at stabilising the concentration 
in certain organs or tissues towards certain desired values (set points)~\cite{SopPatSar14}
and (ii) the drug concentration needs to remain within certain bounds which 
define a \textit{therapeutic window}~\cite{SopPatSarBem15}. 

Another level of classification alludes to 
the mode of administration where we identify the (i) continuous administration 
by, for instance, intravenous continuous infusion, (ii) bolus administration, 
where the medicine is administered at discrete time 
instants~\cite{SopPatSarBem15,RIVADENEIRA2015507} and 
(iii) oral administration where the drug is administered both at discrete times 
and from a discrete set of dosages (e.g., tablets)~\cite{SopSar12}. 

Drug administration is classified according to the way in 
which decisions are made in regard to how often, at what rate and/or what amount 
of drug needs to be administered to the patient. We can identify 
(i) \textit{open-loop} administration policies where the patient follows a 
prescribed dosing scheme without adjusting the dosage and (ii) \textit{closed-loop}
administration where the dosage is adjusted according to the progress of the therapy~\cite{SopPatSar14}.
The latter is suitable mainly for hospitalised patients who are under monitoring and 
where drug concentration measurements are available or the effect of the drug can be
quantified. Such is the case of controlled anaesthesia~\cite{Krieger2014699}. 
However, applications of \textit{closed-loop} policies extend beyond hospitals, 
such as in the case of glucose control in diabetes~\cite{DelFavero2014}. 

Despite the fact that optimisation-based methods are well-established in numerous
scientific disciplines along with their demonstrated advantages over 
other approaches, to date, empirical approaches remain popular%
~\cite{Kannikeswaran20161347,Fukudo2017,deOcenda2016,Savic2017}.

In the following sections we propose decision making approaches for optimal 
drug administration using as a benchmark a fractional two-compartment model 
using, arbitrarily and as a benchmark, the parameters values of amiodarone.
We focus on the methodological framework rather than devising an administration
scheme for a particular compound.
We discuss three important topics in optimal administration
of compounds which are governed by fractional pharmacokinetics.
Next, we formulate the drug administration 
problem as an optimal control problem we prescribe optimal therapeutic courses 
to individuals with known pharmacokinetic parameters.
Lastly, we design optimal administration 
strategies for populations of patients whose pharmacokinetic parameters are 
unknown or inexactly known. We present an advanced closed-loop controlled administration methodology based on model predictive control.

\subsection*{Individualised administration scheduling}%

In this section we present an optimal drug administration scheme based on the two-compartment
fractional model of Eq.~\eqref{eq:model-two-comp}
and assuming that the drug is administered into the central compartment continuously.
We will show that the same model and the same approach can be modified to form the basis 
for bolus administration.

\begin{figure}
 \centering
 \includegraphics[width=0.5\textwidth]{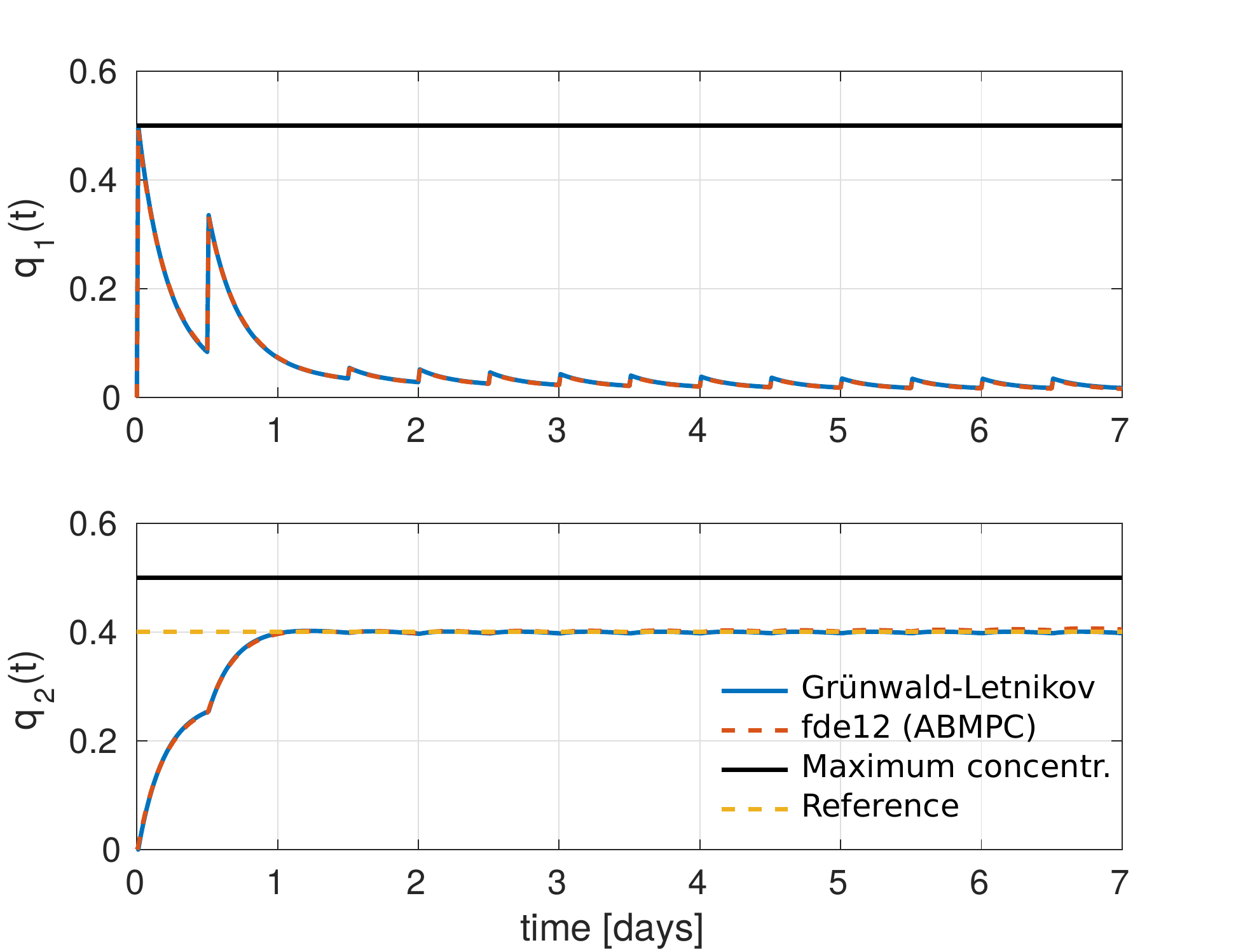}
 \caption{Drug administration scheduling via optimal control. The behaviour of the 
	  controlled system was simulated using the ABMPC method (dashed red line). 
	  The predicted pharmacokinetic profile using the Gr\"unwald-Letnikov approximation 
	  (solid blue line) and \texttt{fde12} are indistinguishable. The set-point on $q_2$ 
	  (dashed orange line) and the maximum allowed concentration (black solid lines) are 
	  shown in the figure.}%
 \label{fig:sp}	  
\end{figure}

In order to state the optimal control problem for optimal drug administration
we first need to discretise the two-compartment model~\eqref{eq:model-two-comp} 
with a (small) sampling time $t_c$ 
\begin{align}
t_c^{-1}(x_{k+1}-x_{k}) = Ax_k + F \,\glDD^{1 - a}_{t_c,\nu} x_k + Bu_k 
\label{eq:asdf-qwerty} 
\end{align}
where $x_k  = \left[ q_1( kt_c) \,\, q_2( kt_c) \right]'$ and $u_k$ is the 
administration rate at the discrete time $k$.
The left hand side of~\eqref{eq:asdf-qwerty} corresponds to the 
forward Euler approximation of the first-order derivative, and we shall refer to 
$t_c=\unit[10^{-2}]{days}$ as the \textit{control sampling time}. 
On the right hand side of~\eqref{eq:asdf-qwerty}, $\glDD^{1 - a}_{t_c,\nu}$ has replaced
the fractional-order operator $_{0}^{\mathrm{c}}\D_t^{1-\alpha}$.
Matrices $A,F$ and $B$ are
\begin{align}
A = \begin{bmatrix}
    -(k_{12} + k_{10}) & 0 \\
      k_{21,f}        & 0 \\
\end{bmatrix},
F = \begin{bmatrix}
     0 & k_{21,f}  \\
     0 &-k_{21,f}  \\
\end{bmatrix},
B = \begin{bmatrix}
     1\\
     0\\
\end{bmatrix}.
\end{align}
The discrete-time dynamic equations of the system can now be stated as 
\begin{align}
    x_{k+1} = x_k 
            + t_c \bigg(
               Ax_k + \frac{F}{t_c^{1-a}} \sum_{j = 0}^{\nu} c_j^{1-a}x_{k-j} 
                    + Bu_k
                 \bigg),
    \label{eq:dynamic_equation}
\end{align}
where $c_j^\alpha = (-1)^j\binom{\alpha}{j} $.
By augmenting the system with past values as 
$\mathbf{x}_k = (x_k,x_{k-1},\dots,x_{k-\nu + 1})$  we can 
rewrite Eq.~\eqref{eq:dynamic_equation} as a finite-dimension linear system
\begin{align}
  \mathbf{x}_{k+1} = \hat{A} \mathbf{x}_k + \hat{B}u_k.
  \label{eq:dynamic_equation_LTI} 
\end{align}
Matrices $\hat{A}$ and $\hat{B}$ are straightforward to derive and are given in~\cite{SopSar16}.
The therapeutic session will last for $N_d = N t_c = \unit[7]{days}$ 
in total, where $N$ is called the \textit{prediction horizon}.
Since it is not realistic to administer the drug to  the patient too frequently, 
we assume that the patient is to receive their treatment every $t_d=\unit[0.5]{days}$.
The administration schedule must ensure that the concentration of 
drug in all compartments never exceeds the minimum toxic concentration 
limits while tracking the prescribed reference value as close as possible.
To this aim we formulate the constrained optimal control problem~\cite{bertsekas-dp}
\begin{subequations}\label{eq:osocp}
 \begin{align}
   \min_{ u_0,\ldots, u_{N_d-1}}J &= \sum_{k=0}^{\nicefrac{N_d}{t_c} + 1} 
      (x_{\mathrm{ref},k} - x_k)' Q (x_{\mathrm{ref},k} - x_k) \\% + \sum_{j=0}^{N_d} \hat{u}_k' R \hat{u}_k \\
 \intertext{subject to the constraints}
  \mathbf{x}_{k+1} &= \hat{A} \mathbf{x}_k + \hat{B} u_j, \text{ for } \ k t_c = j t_d\\
  \mathbf{x}_{k+1} &= \hat{A} \mathbf{x}_k, \text{ otherwise}\\
  0 &\le x_{k} \le 0.5 \\
  0 &\le u_{j} \le 0.5 \\
  \intertext{for $k=0,\ldots, N$; $j=0,\ldots, N_d-1$.} \nonumber
  \end{align}
\end{subequations} 
In the above formulation $x_{\mathrm{ref},k}$ is the desired drug concentration at time $k$ 
and operator $'$ denotes vector transposition.
Any deviation from set point is penalised by a weight matrix $Q$, 
where here we chose $Q= \mathrm{diag}([0\;1])$. 
Note that we are tracking only the second state. 
Our underlying GL model has a relative history of $t_c \nu =\unit[5]{days}$.
Optimal drug concentrations are denoted by $u^\star_k$, for $k = 0,\dots, N_d-1$ and they correspond
to dosages administered intravenously at times $kt_d$. In the optimal control formulation we have implicitly
assumed that $t_d$ is an integer multiple of $t_c$, which is not restrictive since $t_c$ can be chosen arbitrarily.
In Figure~\ref{fig:sp} we present the pharmacokinetic profile of a patient following 
the prescribed optimal administration course.

Finally, the problem in Eq.~\eqref{eq:osocp} is a standard quadratic
problem with simple equality and inequality constraints 
that can be solved at low computational complexity. 
Such problems can be easily formulated in YALMIP~\cite{YALMIP} or ForBES~\cite{forbes} for MATLAB or
CVXPY~\cite{cvxpy} for Python.

The optimal control framework offers great flexibility in using 
arbitrary system dynamics, constraints on the administration rate 
and the drug concentration in the various compartments, and cost 
functions which encode the administration objectives.

The quadratic function which we proposed in Eq.~\eqref{eq:osocp} is certainly not the 
only admissible choice. For instance, other possible choices are 
\begin{enumerate}
 \item $J = \sum_{k=0}^{\nicefrac{N_d}{t_c} + 1} (x_{\mathrm{ref},k} - x_k)' Q (x_{\mathrm{ref},k} - x_k) + \beta u_k^2$, where we
       also penalise the total amount of drug that is administered throughout the prediction horizon,
 \item $J = \sum_{k=0}^{\nicefrac{N_d}{t_c} + 1} \mathrm{dist}^2(x_k, [\underbar{$x$}, \bar{x}])$, 
       where $[\underbar{$x$}, \bar{x}]$ is the therapeutic window and $\mathrm{dist}^2(\cdot, [\underbar{$x$}, \bar{x}])$
       is the squared distance function defined as $\mathrm{dist}^2(x_k, [\underbar{$x$}, \bar{x}]) = \min_{z_k} \{\|x_k-z_k\|^2;\ 
       \underbar{$x$} \leq z_k \leq \bar{x}\}$.
\end{enumerate}

\subsection*{Administration scheduling for populations}%
In the previous section, we assumed that the pharmacokinetic parameters are 
known aiming at an individualised dose regimen. When designing an administration 
schedule for a population of patients --- without the ability to monitor the 
distribution of the drug or the progress of the therapy --- then
$J$ becomes a function of the pharmacokinetic parameters 
(in our case $k_{10}$, $k_{12}$, $k_{21,f}$ and $\alpha$).
That said, $J$ becomes a random quantity which follows an --- often 
unknown or inexactly known --- probability distribution.
In order to formulate an optimal control problem we now need to 
extract a \textit{characteristic value} out of the random quantity $J$.
There are two popular ways to do so leading to different problem statements.

First, we may consider the maximum value of $J$, $\max J$, over all 
possible values of $k_{10}$, $k_{12}$, $k_{21,f}$ and $\alpha$.
This leads to \textit{robust} or \textit{minimax} problem formulations
which consist in solving~\cite{bertsekas-dp}
 \begin{align}\label{eq:minimax-problem}
   \min_{ u_0,\ldots, u_{N_d-1}} \max_{k_{10}, k_{12}, k_{21,f},\alpha} J(u_0,\ldots, u_{N_d-1}, k_{10}, k_{12}, k_{21,f},\alpha),
  \end{align}
subject to the system dynamics and constraints. In Eq.~\eqref{eq:minimax-problem},
the maximum is taken with respect to the worst-case values of $k_{10}$, $k_{12}$, $k_{21,f}$ and $\alpha$,
that is, their minimum and maximum values. Apparently, the minimax approach does not 
make use of any probabilistic or statistical information which is typically available
for the model parameters. As a result, it is likely to be overly conservative
and lead to poor performance.

\begin{figure}[h]
\centering
\includegraphics[width = 0.5\textwidth]{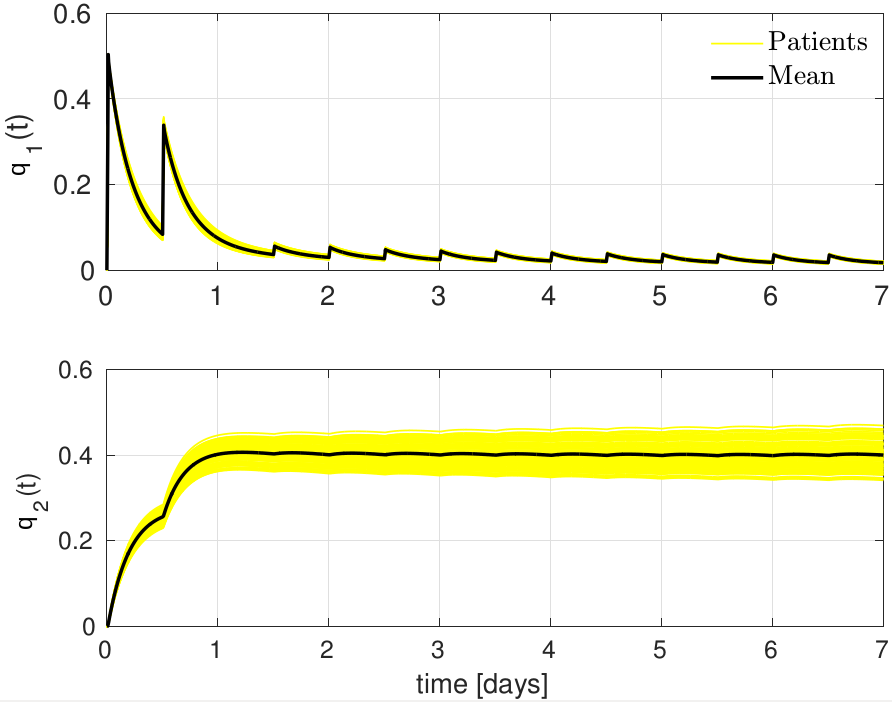}
\caption{Drug administration to a population of $100$ (randomly generated) individuals by stochastic optimal control.
	The yellow thin lines correspond to the predicted individual responses to the administration and the black lines 
	represent the average response.}
\label{fig:sim_multiple_response}
\end{figure}

On the other hand, in the \textit{stochastic} approach we minimise the~\cite{bertsekas-shreve-soc}
expected cost $\mathbb{E}J$ 
 \begin{align}\label{eq:stochastic-problem}
   \min_{ u_0,\ldots, u_{N_d-1}} \mathbb{E}_{k_{10}, k_{12}, k_{21,f},\alpha}\ J(u_0,\ldots, u_{N_d-1}, k_{10}, k_{12}, k_{21,f},\alpha),
  \end{align}
subject to the system dynamics and constraints.
Open-loop stochastic optimal control methodologies have been proposed for the 
optimal design of regimens under uncertainty for classical integer-order pharmacokinetics%
~\cite{Schumitzky83,Lago1992,Baryard1994}, yet, to the best of our knowledge, no studies
have been conducted for the effectiveness of stochastic optimal control for fractional 
pharmacodynamics.

The expectation in Eq.~\eqref{eq:stochastic-problem} can be evaluated 
empirically given a data-set of estimated pharmacokinetic parameters 
$(k_{10}^{(i)}, k_{12}^{(i)}, k_{21,f}^{(i)},\alpha^{(i)})_{i=1}^{L}$
by minimising the sample average of $J$, that is
 \begin{align}\label{eq:stochastic-problem-sampled}
   \min_{ u_0,\ldots, u_{N_d-1}} \tfrac{1}{L}\sum_{i=1}^{L} J(u_0,\ldots, u_{N_d-1}, k_{10}^{(i)}, k_{12}^{(i)}, k_{21,f}^{(i)},\alpha^{(i)}).
  \end{align}
The effect of $L$ on the accuracy of the derived optimal decisions is
addressed in~\cite{Campi2009} where probabilistic bounds are provided.
Simulation results with the stochastic optimal control methodology 
on a population of $100$ patients using $L=20$ are shown in Figure~\ref{fig:sim_multiple_response}.
Note that in both the minimax and the stochastic approach, a common 
therapeutic course $u_0,\ldots, u_{N_d-1}$ is sought for the whole
population of patients.

In stochastic optimal control it is customary to adapt the state 
constraints in a probabilistic context requiring that they be 
satisfied with certain probability, that is
\begin{align}
 \label{eq:chance-constraints}
 \mathrm{P}[x_k \leq x_{\max}] \geq 1-\beta,
\end{align}
where $\beta \in [0,1]$ is a desired level of confidence.
Such constraints are known as \textit{chance constraints} or
\textit{probabilistic constraints}.
Such problems, unless restrictive assumptions are imposed on the 
distribution of $x_k$, lead to nonconvex optimisation problems and 
solutions can be approximated by Monte Carlo sampling methods%
~\cite[Sec.~5.7.1]{shapiro2009lectures} or by means of convex relaxations
known as \textit{ambiguous chance constraints} where~\eqref{eq:chance-constraints}
are replaced by $\mathrm{AV@R}_{\beta}[x_k] \leq x_{\max}$, where 
$\mathrm{AV@R}_{\beta}[x_k]$ is the average value-at-risk of $x_k$ at level $\beta$%
~\cite[Sec.~6.3.4]{shapiro2009lectures}.

The minimax and the expectation operator are the two ``extreme'' choices 
with the former relating to complete lack of probabilistic information and the
latter coming with the assumption of exact knowledge of the underlying 
distribution of pharmacokinetic parameters. Other operators can be 
chosen to account for imperfect knowledge of that distribution, therefore,
bridging the gap between minimax and stochastic optimal control.
Suitable operators for optimal control purposes are the \textit{coherent risk measures} 
which give rise to \textit{risk-averse optimal control} which is the state of 
the art in optimsation under uncertainty~\cite{HerSopBemPat17}.

\subsection*{Model Predictive Control}%
Model predictive control fuses an optimal control open-loop decision making process
with closed-loop control when feedback is available. Following its notable success and 
wide adoption in the industry, MPC has been proposed for several cases of drug administration
for drugs that follow integer-order kinetics such as erythropoietin~\cite{Gaweda2008},
lithium ions~\cite{SopPatSarBem15},
propofol for anaesthesia~\cite{propofol2008} and, most predominantly, insulin 
administration to diabetic patients~\cite{Schaller2016,Hovorka2004,Toffanin2013,Parker1999}.

In MPC, at every time instant, we solve an optimal control problem which produces a \textit{plan}
of future drug dosages over a finite \textit{prediction horizon} by minimising a performance index $J$. 
The future distribution of the drug in the organism is predicted with a dynamical model such as the 
ones based on the Gr\"unwald-Letnikov discretisation presented in the previous sections.
Out of the planned sequence of dosages, the first one is actually administered to the 
patient. At the next time instant, the drug concentration is measured and the same procedure is
repeated in a fashion known as \textit{receding horizon control}~\cite{RawMay09}.

A fractional variant of the classical PID controller, namely a PI$^{\lambda}$D$^\mu$ fractional-order
controller, has been proposed in~\cite{SopSar14} for the controlled administration.
However, the comparative advantage of MPC is that it can inherently account for administration rate, 
drug amount and drug concentration constraints which are of capital importance.
In~\cite{SopNtoSar15,SopSar16}, the truncated Gr\"unwald-Letnikov approximation serves as the basis
to formulate MPC formulations with guaranteed asymptotic stability and satisfaction of the 
constraints at all time instants.  
We also single out impulsive MPC, which was proposed in~\cite{SopPatSarBem15}, which is 
particularly well-suited for applications of bolus administration where the patient is
injected the medication (e.g., intravenously), thus abruptly increasing the drug 
concentration. In such cases, it is not possible to achieve equilibrium, but instead 
the objective becomes to keep the concentration within certain therapeutic limits.
Model predictive control can further be combined with 
stochastic~\cite{PatSopSarBem14,SopHerPatBem17} and risk-averse~\cite{HerSopBemPat17} 
optimal control aiming at a highly robust administration that is resilient to 
the inexact knowledge of the pharmacokinetic parameters of the patient as well as 
potential time varying phenomena (e.g., change in the PK/PD parameter values due to illness,
drug-drug interactions and alternations due to other external influences). 

Furthermore, in feedback control scenarios, state information is often partially available.
For example, when a multi-compartment model is used, concentration measurement from only a 
single compartment is available. As an example, in the case of multi-compartment physiologically based models, we would 
normally only have information from the blood compartment.
In such cases, a \textit{state observer} can be used to produce estimates of concentrations 
or amounts of drug in compartments where we do not have access. 
A state observer is a dynamical system which, using observations of (i) those concentrations
that can be measures and (ii) amount or rate of administered drug at every time instant, 
produces estimates $\hat{q}_i$ of those amounts/concentrations of drug to which we do not 
have access. State observers are designed so that $\hat{q}_i(t) - q_i(t)\to 0$ as $t\to\infty$,
that is, although at the beginning (at $t=0$) we only have an estimate $\hat{q}_i(0)$ of $q_i(0)$,
we shall eventually obtain better results which converge to the actual concentrations.
State observers, such as the Kalman filter, have been successfully employed in the administration 
of compounds following integer-order dynamics~\cite{SopPatSar14} including applications of 
artificial pancreas~\cite{diabetesLQG2007,Wang2014}. 

Observers are, furthermore, employed to filter out measurement errors and modelling errors.
When the system model is inexact, the system dynamics in Eq.~\eqref{eq:dynamic_equation_LTI} becomes
\begin{align}
  \mathbf{x}_{k+1} = \hat{A} \mathbf{x}_k + \hat{B}u_k + d_k, 
\end{align}
where $d_k$ serves as a model mismatch term. We may then assume that $d_k$ follows
itself some dynamics often as trivial as $d_{k+1}=d_k$. We may then formulate the following
linear dynamical system with state $(\mathbf{x}_k,d_k)$
\begin{align}
  \begin{bmatrix}
  \mathbf{x}_{k+1}\\d_{k+1} 
  \end{bmatrix}
  = \begin{bmatrix}
\hat{A} & I\\
0 & I
    \end{bmatrix}
\begin{bmatrix}
      \mathbf{x}_k\\d_k
    \end{bmatrix}
    +
    \begin{bmatrix}
    \hat{B}\\0
    \end{bmatrix}
u_k,
\end{align}
and build a state observer for $(\mathbf{x}_k,d_k)$ jointly.
The estimates $(\hat{\mathbf{x}}_k,\hat{d}_k)$ are then fed to the MPC.
This leads to an MPC formulation known as \textit{offset-free MPC}~\cite{SopPatSar14}
that can control systems with biased estimates of the pharmacokinetic parameters.

Model predictive control is a highly appealing approach because it 
can account for imperfect knowledge of the pharmacokinetic parameters
of the patient, measurement noise, partial availability of measurements
and constraints while it decides the amount of administered drug by 
optimising a cost function thus leading to unparalleled performance. 
Moreover, the tuning of MPC is more straightforward compared to 
other control approaches such as PID or fractional PID. In MPC the main 
tuning knob is the cost function $J$ of the involved optimal control 
problem which reflects and encodes exactly the control objectives.

\section*{Conclusions}
Fractional kinetics offers an elegant description of anomalous kinetics, i.e., non-exponential terminal 
phases, the presence of which has been acknowledged in pharmaceutical literature extensively. The 
approach offers simplicity and a valid scientific basis, since it has been applied in problems of 
diffusion in physics and biology. It introduces the Mittag-Leffler function which describes power law 
behaved data well, in all time scales, unlike the empirical power-laws which describe the data only for 
large times. Despite the mathematical difficulties, fractional pharmacokinetics offer undoubtedly a powerful 
and indispensable approach for the toolbox of the pharmaceutical scientist. 

Solutions of fractional-order systems involve Mittag-Leffler-type 
functions, or other special functions whose numerical evaluation is 
nontrivial. Several approximation techniques have been proposed as 
we outlined above for fractional systems with 
different levels of accuracy, parsimony and suitability for optimal administration scheduling 
or control. 
In particular, those based on discrete time-domain approximations of the system 
dynamics with bounded approximation error, such as the truncated Gr\"unwald-Letnikov approximation,
are most suitable for optimal control applications.
There nowadays exist software that implement most of the algorithms that are available
in the literature and facilitate their practical use.
Further research effort needs to be dedicated in deriving error bounds in the 
time domain for available approximation methods that will allow their adoption 
in optimal control formulations.

The optimal control framework is suitable for the design of administration courses both to 
individuals as well as to populations where the intra-patient variability of pharmacokinetic
parameters needs to be taken into account. In fact, stochastic optimal control offers a 
data-driven decision making solution that enables us to go from sample data to 
administration schemes for populations.
Optimal control offers great flexibility in 
formulating optimal drug dosing problems and different problem structures arise naturally
for different modes of administration (continuous, bolus intravenous, oral and more).
At the same time, further research is necessary to make realistic assumptions and translate 
them into meaningful optimisation problems.

MPC methodologies are becoming popular for controlled drug 
administration. Yet, their potential for fractional-order pharmacokinetics
and their related properties, especially in regard to stochastic systems and 
the characterisation of invariant sets, needs to be investigated.

Overall, at the intersection of fractional systems theory, pharmacokinetics, 
numerical analysis, optimal control and model predictive control, spring numerous 
research questions which are addressed in the context of mathematical 
pharmacology.

\bibliographystyle{spbasic}
\bibliography{bibliography}

\end{document}